\documentclass[11pt]{article}
\usepackage{graphicx}
\usepackage{amsmath,latexsym,amsbsy,amssymb, color}
\usepackage{graphicx}
\font\bigbf=cmbx10 at 16pt

\def\ds{\displaystyle}
\def\forall{\hbox{for all}~}
\def\L{{\bf L}}

\def\ve{\varepsilon}

\def\J{{\cal J}}
\def\I{{\cal I}}

\def\R{{\mathbb R} }

\def\vp{\varphi}

\def\prob{\hbox{Prob.}}

\def\vs{\vskip 2em}
\def\vsk{\vskip 4em}
\def\v{\vskip 1em}

\def\S{{\cal S}}
\def\C{{\cal C}}

\def\W{{\cal W}}

\def\ov{\overline}
\def\wto{\rightharpoonup}
\def\Tilde{\widetilde}
\def\Hat{\widehat}
\def\bega{\begin{array}}
\def\enda{\end{array}}
\def\begi{\begin{itemize}}
\def\endi{\end{itemize}}
\def\meas{\hbox{meas}}

\def\P{\hbox{Prob.}}
\def\bel{\begin{equation}\label}
\def\eeq{\end{equation}}
\def\sqr#1#2{\vbox{\hrule height .#2pt
\hbox{\vrule width .#2pt height #1pt \kern #1pt
\vrule width .#2pt}\hrule height .#2pt }}
\def\square{\sqr74}
\def\endproof{\hphantom{MM}\hfill\llap{$\square$}\goodbreak}

\newtheorem{theorem}{Theorem}[section]

\newtheorem{lemma}{Lemma}[section]

\newtheorem{definition}{Definition}[section]

\begin{document}

\title{\bigbf Diffusion Approximations of 
Markovian Solutions to Discontinuous ODEs}
\vs
\author{Alberto Bressan$^{(1)}$, Marco Mazzola$^{(2)}$, and Khai T. Nguyen$^{(3)}$\\ 
\\
 {\small $^{(1)}$ Department of Mathematics, Penn State University, }\\  
 {\small $^{(2)}$ IMJ-PRG, CNRS, Sorbonne Universit\'e,}\\
  {\small $^{(3)}$ Department of Mathematics, North Carolina State University. }\\ 
 \\  {\small e-mails: ~axb62@psu.edu, ~marco.mazzola@imj-prg.fr,  ~ khai@math.ncsu.edu}
 }
\maketitle

\begin{abstract}  In a companion paper, the authors have characterized
all deterministic semigroups, and all Markov semigroups, 
whose trajectories are Carath\'eodory
solutions to a given ODE  $\dot x = f(x)$, with $f$ possibly discontinuous.
The present paper establishes two approximation results.  
Namely, every deterministic
semigroup can be obtained as the pointwise limit of the flows 
generated by a sequence of ODEs $\dot x=f_n(x)$ with smooth right hand sides.
Moreover, every Markov semigroup can be obtained as
limit of a sequence of diffusion processes with smooth drifts and with diffusion coefficients approaching zero.
\end{abstract}

\vsk

\section{Introduction}
\label{s:1}
\setcounter{equation}{0}

Consider the Cauchy problem for a scalar  ODE with possibly
discontinuous  right hand side:
\bel{ode}
\dot x~=~f(x),\eeq
\bel{ic}x(0)~=~x_0\,.\eeq
By definition, a map $t\mapsto x(t)$ is a {\em Carath\'eodory solution} 
of (\ref{ode})-(\ref{ic})  if
\bel{CS}
x(t)~=~x_0+\int_0^tf(x(s))\, ds\qquad\qquad \forall t\geq 0\,.
\eeq
When the function $f$ is not Lipschitz continuous, it is well known that 
this Cauchy problem can admit multiple solutions.
Because of this non-uniqueness, 
in \cite{E, FL} it was proposed to study ``generalized flows", described by a probability
measure on the set of all Carath\'eodory solutions. In this direction,  in the companion paper \cite{BMN} the authors
have characterized all deterministic semigroups, and all  Markov semigroups,
whose trajectories are solutions to the ODE in (\ref{ode})-(\ref{ic}).

Aim of the present paper is to show that all of these semigroups
can be obtained as limits of smooth approximations. To explain these results 
more precisely, we recall

\begin{definition}\label{d:DS} A {\bf deterministic semigroup} compatible with the
ODE (\ref{ode}) is a map $S:\R\times \R_+\mapsto \R$,
with the properties
\begi\item[(i)] $S_t(S_s(x_0))~=~S_{t+s}(x_0)$, ~$S_0(x_0)~=~ x_0$.
\item[(ii)] For each $x_0\in \R$, the map $t\mapsto S_t x_0$
is a solution to the Cauchy problem (\ref{ode})-(\ref{ic}).
\endi
\end{definition}
Notice that here we do not require any continuity w.r.t.~the initial point $x_0$.
Of course, if $f$ is Lipschitz continuous, then the solution 
$t\mapsto x(t) = S_t x_0$  of (\ref{ode})-(\ref{ic}) is unique, and this
uniquely determines the semigroup.

Throughout this paper, as in \cite{BMN} we shall assume:
\begi
\item[{\bf (A1)}] {\it The function $f:\R\mapsto\R$ is  bounded and regulated.
Namely, $f$  admits left and right limits $f(x-)$, $f(x+)$ at every point $x$, while $M\doteq 
\sup_x |f(x)|<+\infty$.}
\item[{\bf (A2)}] {\it If $y$ is  a point  where
either $f(y-)\cdot f(y+)=0$ or else $f(y-)> 0> f(y+)$, then $f(y)=0$.}
\endi
As in \cite{Binding}, the ``no jam" assumption {\bf (A2)} guarantees that the Cauchy problem (\ref{ode})-(\ref{ic}) has at least one solution. 
Indeed, in the present setting every solution in the sense of Filippov \cite{Fil} is a Carath\'eodory solution as well.
Under the above assumptions, 
the set of discontinuities
\bel{Df}
D_f~=~\Big\{x\in\R\,;~f(x-)\neq f(x+)~~\mathrm{or}~~f(x-)=f(x+)\neq f(x)\Big\}
\eeq
is at most countable. Moreover,  the set
of zeroes
\bel{zset}f^{-1}(0)~=~\bigl\{x\in\R\,;~f(x)=0\bigr\}\eeq
is closed.  
Given a regulated function $f:\R\mapsto\R$, we can complete its graph by
adding a vertical segment at each point $x$ where $f$ has a jump.
This yields a multifunction with closed graph and 
compact, convex values:
\bel{F}F(x)~\doteq~\mathrm{co}\bigl\{f(x), \, f(x+),\, f(x-)\bigr\},\eeq
where $``\mathrm{co}"$ denotes the convex closure.
Throughout the sequel, 
$B(V,r)=\{x\,;~d(x,V)<r\}$  denotes the open neighborhood of radius $r$ 
around the set $V$.

\begin{definition}\label{d:gc} Let $f$ be a regulated function and 
let $F$ be the corresponding multifunction in (\ref{F}). 
We say that a sequence of smooth functions
$f_n$ converges  $f$ in the sense of the graph if, 
for every given $\ve>0$ and $N>0$, one has
\bel{graph}\hbox{\rm Graph}\Big(f_n\Big|_{[-N,N]}\Big)~\subset~B\bigl(\hbox{\rm Graph}(F)\,,~\ve\bigr)
\eeq
for all $n$ sufficiently large. \end{definition}
In other words, for all $n$ sufficiently large,
every point $(x, f_n(x))$, with $|x|\leq N$, 
is contained in an $\ve$-neighborhood
of the graph of $F$.  We recall that this approach has been used in the literature, 
in the analysis of upper semicontinuous multifunctions with convex values \cite{AC, Ce}.
\medskip

Our first main result shows that every deterministic semigroup $S$ can be 
approximated by the semigroups $S^n$ generated  by a sequence of  smooth ODEs.
\bel{CPn}\dot x ~=~ f_n(x),\qquad\qquad x(0)~=~ x_0\,.\eeq

\begin{theorem}\label{t:1} Let $f$ satisfy the assumptions {\bf (A1)-(A2)}.
 Let $S:\R\times \R_+\mapsto\R$ 
be a semigroup compatible with the ODE (\ref{ode}).   Then there exists
a sequence of smooth functions $(f_n)_{n\geq 1}$, converging to $f$ 
in the sense of the graph, such that the following holds. For each $x_0\in\R$, calling $t\mapsto S^n_t x_0$ the  unique solution to the 
Cauchy problem (\ref{CPn}),
one has the convergence
\bel{limn}\lim_{n\to\infty}\,\sup_{t\geq 0}~
\bigl|S^n_t x_0-S_t x_0\bigr|~=~0.\eeq
\end{theorem}
 We remark that,
since the map 
$x_0\mapsto S_t x_0$ can be discontinuous, one can only achieve the
pointwise convergence for each fixed initial datum $x_0\in\R$.
On the other hand, we can achieve uniform convergence w.r.t. the time variable 
$t\in [0,+\infty[\,$.

The second part of the paper is concerned with Markov semigroups.
We consider a family of transition kernels
$P_t(x_0,A)$, denoting the probability that a solution starting at $x_0$
reaches a point in the Borel set $A\subset\R$ at time $t$.
It is assumed that these transition kernels satisfy the Chapman-Kolmogorov equation
\bel{CK}
P_{t+s}(x_0,A)~=~\int P_t(z,A)\, P_s(x_0, dz).
\eeq
\begin{definition}\label{d:13}
We say that a family of transition kernels $P_t(\cdot,\cdot)$ defines  a Markov semigroup 
 compatible with the ODE (\ref{ode}) if the following holds.
There exists a probability space 
${\mathcal W}$ and a Markov process $X(t, x_0,\omega)$, $\omega\in\W$,  such that, for every $x_0\in\R$ and 
every Borel set $A\subset\R$,
\bel{MP}
\P\Big\{ X(t,x_0,\omega)\in A\Big\}~=~P_t(x_0,A),\eeq
and moreover all sample paths $t\mapsto X(t,x_0,\omega)$, $\omega\in 
{\mathcal W}$, are Carath\'eodory solutions to Cauchy problem (\ref{ode})-(\ref{ic}).   
\end{definition}

Our second main result shows that 
every one of these Markov semigroups 
can be obtained as a limit of a sequence of diffusion processes, of the form
\bel{diffu} 
dX~=~f_n(X)\, dt + \sigma_n dW.\eeq
Here $(f_n)_{n\geq 1} $ is a sequence of smooth functions, 
converging to $f$ in the sense of the graph, while the diffusion coefficients 
$\sigma_n$ decrease to zero, as $n\to\infty$. 
As usual, by $W$ we denote the standard one-dimensional
Wiener process.

It is well known that the transition probability kernel $P^{(n)}$
for (\ref{diffu}) can be obtained as follows. Let $\Gamma_n(t, x;x_0)$ be the fundamental solution to the 
linear parabolic equation
\bel{pe} v_t + \bigl(f_n(x)v\bigr)_x~=~{\sigma^2_n\over 2}\, v_{xx}\,,\eeq
where the initial data is a unit mass at the point $x_0$, namely
\bel{uxo}\lim_{t\to 0+} \Gamma_n(t,x;x_0)~=~\delta_{x_0}\eeq
in distributional sense.
Then, for every $t>0$ and every Borel set $A\subset\R$, 
\bel{Pnt}
P^{(n)}_t(x_0, A)~=~\int_A \Gamma_n(t,y; x_0)\, dy.\eeq
The next theorem shows that, by a suitable choice of the drifts $f_n$,
one can achieve the convergence $P^{(n)}\to P$ in distribution.
\begin{theorem}\label{t:2} Let $f$ satisfy the assumptions {\bf (A1)-(A2)}.
Let $P_t(x_0,A)$ be a family of transition kernels, defining a Markov 
semigroup compatible with the ODE (\ref{ode}). 
Then there exists a sequence of diffusion processes of the
form (\ref{diffu}) with $f_n\in \C^\infty$, such that
\begi
\item[(i)] $f_n\to f$ in the sense of the graph,
\item[(ii)] $\sigma_n\to 0$,
\item[(iii)] for every $x_0\in \R$ and $t>0$, the transition kernels converge
in distribution.  Namely, for every bounded continuous function $\vp\in \C^0(\R)$, 
one has the convergence of the expected values
\bel{EP}
\lim_{n\to\infty} 
\int \vp(x)P^{(n)}_t(x_0,dx)~=~\int \vp(x) P_t(x_0, dx).
\eeq
\endi
\end{theorem} 

Here the heart of the matter is to construct a sequence of piecewise constant drifts $f_n\in \L^\infty$ 
and diffusion coefficients $\sigma_n>0$ which satisfy (i)--(iii).    The additional 
property $f_n\in \C^\infty$ is then achieved by suitable mollifications.
We again remark that,
since in general the limit semigroup has no continuous dependence w.r.t.~the initial point $x_0$,
the convergence $P^{(n)}_t(x_0,\cdot)\to P_t(x_0, \cdot)$ can only be attained 
in a pointwise sense, for each initial point $x_0$.

The remainder of the paper is organized as follows. Section~\ref{s:2} reviews 
 some earlier results on discontinuous ODEs, while Section~\ref{s:3}
contains  a proof of Theorem~\ref{t:1}.
Finally, Theorem~\ref{t:2} is proved in Sections~\ref{s:4} to \ref{s:6}.

To explain the main ideas involved in these proofs we recall that, as proved in \cite{BMN},
to single out a unique deterministic semigroup compatible with (\ref{ode}),
three additional ingredients are needed:
\begi
\item[{\bf (Q1)}]
A continuum (i.e., atomless) positive measure $\mu$ supported on 
the set of zeroes (\ref{zset}) of $f$.
Strictly increasing  trajectories $t\mapsto x(t)= S_t(x_0)$ of the semigroup are then implicitly defined by the identity
$$t~=~\int_{x_0}^{x(t)} {dy\over f(y)}+ \mu\bigl([x_0, x(t)]\bigr),$$
while a similar formula holds for decreasing ones.
Notice that, if $f^{-1}(0)$ is countable, then necessarily $\mu=0$.
However, in \cite{BMN} an example was constructed where $f^{-1}(0)\subset [0,1]$ is 
the Cantor set, while $f$ is H\"older continuous and strictly positive at all other points.    In this case, different choices of the measure $\mu$ lead to infinitely many different semigroups, all compatible with the ODE.

\item[{\bf (Q2)}] A countable set of points $\S\subseteq f^{-1}(0)$, where the dynamics
is forced to stop.   Among the (possibly many) solutions of (\ref{ode})
starting from  $x_0\in \S$, this means that we are selecting the stationary one:
$S_t (x_0)=x_0$.
\item[{\bf (Q3)}] A map $\Phi:\Omega^*\mapsto \{-1,1\}$, defined
on a set $\Omega^*\subset\R$ of isolated points from where
both an increasing and a decreasing solution of (\ref{ode}) can originate.
For $x_0\in \Omega^*$, setting $\Phi(x_0)= 1$ selects the increasing
solution, while $\Phi(x_0)=-1$ selects the decreasing one.
\endi

To prove Theorem~\ref{t:1}, 
 we first identify
maximal open intervals $J_k=\,]x_{k-1}, x_k[$ where
trajectories
of the semigroup $S$ are strictly increasing, or strictly decreasing.
Then, to achieve the convergence (\ref{limn}), we 
construct a sequence of smooth functions $f_n$ with the following property.
For each interval $J_k$ where the dynamics is increasing and every two points
$a,b\in J_k$, if  $b=S_\tau a$ for some $\tau=\tau(a,b)$, then
\bel{TC}
\lim_{n\to\infty} \int_a^b {1\over f_n(x)} \, dx~=~\tau.\eeq
Notice that (\ref{TC}) yields the convergence of the times needed for trajectories
to move from $a$ to $b$. We choose the $f_n$ so that the same property
also holds on intervals  where the dynamics is strictly decreasing.

Toward a proof of Theorem~\ref{t:2} we recall that, still by the results in \cite{BMN},
a general  Markov semigroup compatible with (\ref{ode})
is determined  by adding two more items to the 
list {\bf (Q1)--(Q3)}. Namely:
\begi
\item[{\bf (Q4)}]  A countable set $\S^*\subset f^{-1}(0)$ and a
map $\Lambda: \S^*\mapsto [0, +\infty]$, describing the random
waiting time of a trajectory which reaches a point $x_k\in\S^*$.
More precisely,  a solution initially at $x_k\in \S^*$ remains at 
$x_k$ for a random time $T_k\geq 0$, then starts moving.
All these random waiting times are mutually independent, with Poisson distribution:
\bel{Tk}\prob\{T_k > s\} ~=~e^{-\lambda_k s}\qquad\mathrm{with}\qquad \lambda_k~=~\Lambda(x_k). \eeq
\item[{\bf (Q5)}] A map $\Theta:\Omega^*\mapsto [0,1]$, defined 
on the countable set $\Omega^*\subset\R$ of points from which
both an increasing and a decreasing solution of (\ref{ode}) can originate.
For $x_k\in \Omega^*$, the value $\theta_k =\Theta(x_k)$ gives the probability that,
when the solution starting from $x_k$ begins to move, it will be increasing.
Of course, $1-\theta_k$ is then the probability that the solution will be decreasing.
\endi

\begin{figure}[ht]
\centerline{\hbox{\includegraphics[width=16cm]{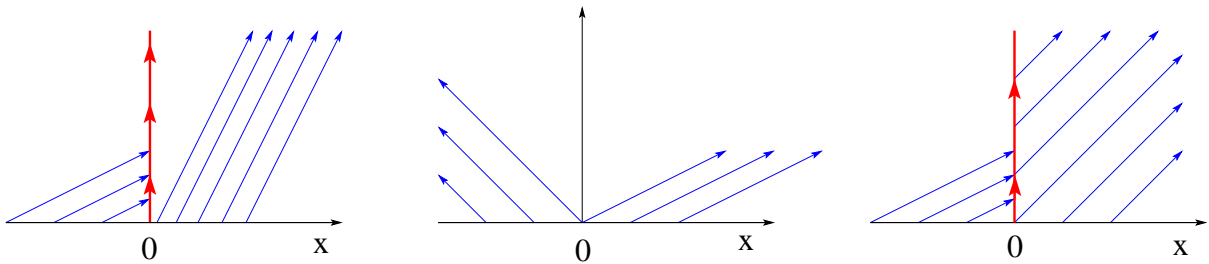}}}
\caption{\small  Left: Case 1, where all trajectories stop at the origin.
Center: Case 2, where trajectories starting at the origin move to the right or to the left with
probabilities $\theta$ and $1-\theta$.   Right: Case 3, where trajectories reaching the origin
wait for a random time, then start moving again to the right. }
\label{f:ode18}
\end{figure}

After various approximations and reductions, discussed in Sections~\ref{s:4} and 
\ref{s:5}, we are led to study
a discontinuous ODE having the basic form
\bel{bode}
\dot x~=~f(x)~=~\left\{\bega{cl} a\quad &\hbox{if}~~x< 0,\cr
0\quad &\hbox{if}~~x= 0,\cr
b\quad &\hbox{if}~~x>0,\enda\right.\eeq
for some $a,b\not= 0$.
For a Markov semigroup compatible with (\ref{bode}), three main cases must be considered, see Fig.~\ref{f:ode18}.

\v
CASE 1: $a,b>0$ and all trajectories reaching  the origin  remain at the origin forever after.
\v

CASE 2: $a<0<b$ and all trajectories starting at the origin  move to the right with probability 
$\theta\in [0,1]$ and to the left with probability $1-\theta$. 
\v
CASE 3: $a,b>0$ and all trajectories reaching the origin wait for a random time $T$
with Poisson distribution
\bel{poisson}\prob\{T > s\} ~=~e^{-\lambda s},\eeq
then start moving to the right with speed $b$.
\v
For the above three cases, 
sequences of piecewise constant drifts $g_n$ achieving the desired convergence 
are constructed in (\ref{gbx}), (\ref{gn*}), and (\ref{gndef}), respectively.

The approximation of a random waiting time (\ref{poisson})   
requires a more careful analysis.   This is based on the 
preliminary construction of a family of eigenfunctions $u_n$ of a corresponding rescaled problem on the interval $[0,1]$, with eigenvalues $\lambda_n\to \lambda$.  
To help the reader, the distribution function for  the limit Markov semigroup 
is shown in Fig.~\ref{f:ode20}. In addition, a lower and and an upper solution to the 
parabolic equation describing the approximating diffusion process 
are shown in Figures~\ref{f:ode21}
and \ref{f:ode22}, respectively.

The convergence of a sequence of stochastic processes to a diffusion process
is a classical topic in probability theory.  See for example \cite{GS, KT} or
Chapter~11 in \cite{SV}.  The present result goes in the converse direction,
approximating a Markov process with
smooth diffusions.  In a way, this is in the same spirit as \cite{CCS}. Indeed, we show
that the existence of diffusion approximations is not a selection
principle for any of the Markov semigroups compatible with (\ref{ode}).
A survey of results on solutions to generalized ODEs can be found in
\cite{BFP}.

\section{Review of scalar discontinuous ODEs}
\setcounter{equation}{0}
\label{s:2}
As a preliminary, we collect here
 some basic results on the existence and properties of solutions to 
a  Cauchy problem with possibly discontinuous right hand side.
In the following theorem, 
the first two statements are proved in somewhat greater generality in \cite{Binding}.
On the other hand,  the closure of the solution set strongly relies on the assumption
{\bf (A1)} that the function $f$ is regulated.  See \cite{BMN} for a proof.
\v
\begin{theorem}\label{t:21} Consider the Cauchy problem (\ref{ode})-(\ref{ic}),
assuming that {\bf (A1)-(A2)} hold. Then 
\begi
\item[(i)]  For every $x_0\in\R$ there exists at least one Carath\'eodory  solution, defined for all times $t\geq 0$. 
\item[(ii)]
Every solution is monotone (either increasing or decreasing).
\item[(iii)] If $x_n:[0,\tau]\mapsto \R$ is a sequence of solutions of (\ref{ode}), and the pointwise convergence $x_n(t)\to x(t)$ holds for every $t\in [0,\tau]$, 
then $x(\cdot)$ is a solution as well.
\endi
\end{theorem}

The general form of a semigroup $S$ compatible with the ODE (\ref{ode})
was described in \cite{BMN}. 
To uniquely determine one of these semigroups, together with the ODE
one needs to  
assign  a positive, atomless measure $\mu$ supported on $f^{-1}(0)$, 
a countable set $\S\subseteq f^{-1}(0)$, and a map $\Phi:\Omega^*\mapsto 
\{-1, 1\}$, as described at  {\bf (Q1)--(Q3)} in the Introduction.  

We briefly review the construction in \cite{BMN}. Given $f$, consider the 
regulated functions
\bel{fpm}
f^+(x)~\doteq~\max\{ 0, f(x)\}, \qquad \quad f^-(x)~\doteq~\min\{0, f(x)\}.\eeq

\begin{definition}\label{d:increase} We say that an open interval $]a,b[$ is a 
{\bf domain of increase}
if
\bel{din}
]a,b[~\cap~\S~=~\emptyset,\qquad \mu([c,d])+\int_c^d {dx\over f^+(x)}~<~+\infty
\qquad\forall [c,d]\subset\,]a,b[\,.\eeq
Similarly, we say that  $]a,b[$ is a 
{\bf domain of decrease}
if
\bel{dde}
]a,b[~\cap~\S~=~\emptyset,\qquad \mu([c,d])-\int_c^d {dx\over f^-(x)}~<~+\infty
\qquad\forall [c,d]\subset\,]a,b[\,.\eeq
\end{definition}

If $]a,b[$ and $]a', b'[$ are two intervals of increase having non-empty
intersection, then their union $]a,b[\,\cup\,]a', b'[$ is also an interval of increase.
We can thus identify countably many, disjoint maximal intervals of increase
$]\alpha_i, \beta_i[$, $i\in \I^+$. 
Similarly, we can identify countably many disjoint maximal intervals of decrease
$]\gamma_i, \delta_i[$, $i\in \I^-$.

It remains to analyze what happens at the endpoints of these intervals.
\begi
\item
If $\alpha_i\notin\S$ and, for some $\ve>0$
\bel{ep1} \mu([\alpha_i, \,\alpha_i+\ve])+\int_{\alpha_i}^{\alpha_i+\ve} {dx\over f^+(x)}~<~+\infty,\eeq
 we then consider 
the half-open interval $I^+_i\doteq [\alpha_i, \beta_i[\,$.
Otherwise, we let $I^+_i$ be an open interval:
$I^+_i\doteq \,]\alpha_i, \beta_i[\,$.
\item
If $\delta_i\notin\S$ and, for some $\ve>0$
\bel{ep2} \mu([\delta_i-\ve, \,\delta_i])-\int_{\delta_i-\ve}^{\delta_i} {dx\over f^-(x)}~<~+\infty,\eeq
 we then consider 
the half-open interval $I^-_i\doteq\, ]\gamma_i, \delta_i]\,$.
Otherwise, we let $I^-_i$ be an open interval:
$I^-_i\doteq \,]\gamma_i, \delta_i[\,$.
\endi

For each $i\in\I^+$, we now describe the increasing dynamics on the 
intervals $I^+_i$.   Given $x_0\in I_i^+$, we consider the time
$$\tau^+(x_0)~\doteq~\mu([x_0, \beta_i]) + \int_{x_0}^{\beta_i} {dy\over f^+(y)}~\in~]
0,\,+\infty].$$
We then set
\bel{S+x}
S^+_t(x_0)~\doteq~ x(t),\eeq
where $x(t)$ is implicitly defined by 
\bel{St+}\left\{\bega{cl}\ds
\mu([x_0, x(t)]) + \int_{x_0}^{x(t)} {dy\over f^+(y)}~=~t \qquad 
&\hbox{if}~t<\tau^+(x_0),\\[4mm]
x(t)~=~\beta_i\qquad&\hbox{if}~t\geq\tau^+(x_0).\enda\right.
\eeq
The construction of the decreasing dynamics on the intervals $I_i^-$ is entirely similar.
Given $x_0\in I_i^-$, we consider the time
$$\tau^-(x_0)~\doteq~\mu([\gamma_i, x_0]) - \int_{\gamma_i}^{x_0} {dy\over f^+(y)}~\in~]
0,\,+\infty].$$
We then set
\bel{S-x}
S^-_t(x_0)~\doteq~ x(t),\eeq
where $x(t)$ is implicitly defined by 
\bel{St-}\left\{\bega{cl}\ds
\mu([x(t), x_0]) - \int^{x_0}_{x(t)} {dy\over f^-(y)}~=~t \qquad 
&\hbox{if}~t<\tau^-(x_0),\\[4mm]
x(t)~=~\beta_i\qquad&\hbox{if}~t\geq\tau^-(x_0).\enda\right.
\eeq
We can now combine together the above solutions, and define the semigroup
$S$ on the whole real line, as follows.
\bel{sd1}
S_t(x_0)~=~\left\{\bega{cl} x_0\qquad &\hbox{if}\quad x_0\notin\left(\bigcup_i I_i^+\right)
\cup\left(\bigcup_i I_i^-\right),\\[4mm]
S_t^+(x_0)\qquad &\hbox{if}\quad x_0\in\left(\bigcup_i I_i^+\right)\setminus
\left(\bigcup_i I_i^-\right),\\[4mm]
S_t^-(x_0)\qquad &\hbox{if}\quad x_0\in\left(\bigcup_i I_i^-\right)\setminus
\left(\bigcup_i I_i^+\right).\enda\right.\eeq
To complete the definition, it remains to define $S_t(x_0)$ in the case
$x_0\in I_i^+\cap
 I_j^-$, for some $i\in \I^+$, $j\in \I^-$.   Notice that this can happen 
only if 
$$x_0~=~\alpha_i ~=~ \delta_j\,,$$
where $I_i^+=[\alpha_i, \beta_i[\,$ and $I_j^-= \,]\gamma_j, \delta_j]$
are half-open intervals where the dynamics is increasing and decreasing, respectively.
By our definitions, this implies $x_0\in \Omega^*$.
Recalling {\bf (Q3)}, we thus define
\bel{sd2}
S_t(x_0)~=~\left\{ \bega{rl} S_t^+(x_0)\qquad &\hbox{if}\quad \Phi(x_0) = ~1,\\[3mm]
 S_t^-(x_0)\qquad &\hbox{if}\quad \Phi(x_0) = -1.
 \enda\right.\eeq
\v
The first main result in \cite{BMN} shows that every semigroup compatible with the
ODE (\ref{ode}) has the above form.

\begin{theorem}\label{t:22}  Let $f:\R\mapsto\R$ satisfy the assumptions {\bf (A1)-(A2)}. The following statements are equivalent.
\begi
\item[(i)] The map $S:\R\times \R_+\to\R$ is a deterministic semigroup compatible with the ODE (\ref{ode}).
\item [(ii)] There exist a positive atomless measure $\mu$ supported on the 
set $f^{-1}(0)$, a countable set of points $\S\subseteq f^{-1}(0)$, and a map $\Phi:\Omega^*\mapsto\{-1,1\}$
as in {\bf (Q1)--(Q3)} such that $S$ coincides with the corresponding 
semigroup constructed at (\ref{S+x})--(\ref{sd2}).
\endi
\end{theorem}
Next, we consider a Markov semigroup whose sample paths
satisfy the ODE (\ref{ode}) with probability one.
To uniquely determine such a Markov process, in addition to {\bf (Q1)-(Q2)}
one needs  to assign:
\begi
\item A countable set $\S^*\subseteq f^{-1}(0)$ and, for each $y_j\in \S^*$, a 
number $\lambda_j>0$ characterizing
the Poisson
waiting time, when  trajectories reach the point $y_j\in \S^*$, as in {\bf (Q4)}.
\item A map $\Theta:\Omega^*\mapsto [0,1]$ determining the probability of moving 
upwards or downwards, from an initial point $z_k\in \Omega^*$, as in 
{\bf (Q5)}.
\endi

We start by selecting an underlying probability space $\W$ such that,
as $\omega\in \W$, the countably many random variables $Y_j(\omega)\in \R_+$ and $Z_k(\omega)\in \{-1,1\}$ are independent, with distributions 
\bel{pxj}
\P\{ Y_j>s\}~=~e^{-\lambda_j s}\quad\forall s>0,\quad
\qquad \quad \P\{ Z_k=1\} ~=~\Theta(z_k).\eeq

Next, we recall that, given the measure $\mu$ and the countable set $\S$,
as above one can uniquely determine the sets 
\bel{om+-}\Omega^+~=~\bigcup_{i\in\I^+} I_i^+\,,\qquad\qquad 
\Omega^-~=~\bigcup_{i\in\I^-} I_i^-\eeq
consisting of countable unions of disjoint intervals where the dynamics can increase or decrease, respectively.
For $x_0\in \Omega^+$, a trajectory $t\mapsto S_t^+(x_0)$ was defined at 
(\ref{S+x})-(\ref{St+}), while for 
$x_0\in \Omega^-$, a trajectory $t\mapsto S_t^-(x_0)$ was defined at 
(\ref{S-x})-(\ref{St-}).   

To construct a Markov process with sample paths $t\mapsto X(t,x_0,\omega)$
satisfying (\ref{ode})-(\ref{ic}), 
for  every
$t>0$ we need to define the transition probabilities
$$
P_t(x_0, A)~=~\P\Big\{X(t,x_0,\omega)\in A\Big\}
,$$
for any initial point $x_0\in \R$ and any Borel set $A\subset\R$.
Recalling our construction of a deterministic flow 
at (\ref{sd1})-(\ref{sd2}), this can be done as follows.
\begi
\item[(i)] If $x_0\notin\Omega^+\cup\Omega^-$, then 
$X(t, x_0, \omega) = x_0$ for every $t\geq 0$ and $\omega\in \Omega$.
Hence 
\bel{TP1} P_t(x_0, \{x_0\})~ = ~1\qquad \forall ~t\geq 0.\eeq
In other words, all  trajectories starting at $x_0$ remain constant.
\item[(ii)] If $x_0\in  \Omega^+\setminus\Omega^-$, a random trajectory starting at 
$x_0$ will have the form
\bel{X+}
X(t, x_0, \omega)~=~S^+_{T^+(t,x_0,\omega)} (x_0)\,,\eeq
where the time $T^+$ along the trajectory is a random variable with distribution 
\bel{T+} {\P\bigl\{ T^+(t, x_0, \omega)<s\bigr\}}~=~\P\left\{ s+
\sum_{y_j\in S^*\cap [ x_0, \, S^+_s(x_0)]}  Y_j (\omega)~>~t\right\}.\eeq
Note that (\ref{T+}) accounts for the (possibly countably many) waiting times
when the trajectory crosses one of the points $y_j\in \S^*$.
The transition probabilities are thus given by
\bel{TP2}\bega{l}
P_t\bigl( x_0,\, \,]-\infty, x_0[\bigr) ~=~0,\\[3mm]
\ds P_t\bigl( x_0,\, \,[x_0, S_s^+(x_0)]\bigr)~=~\P\left\{ s+\sum_{y_j\in S^*\cap [ x_0, \, S^+_s(x_0)]}  Y_j (\omega)~\geq ~t\right\}.
\enda \eeq

\item[(iii)] Similarly, for an initial state  $x_0\in  \Omega^-\setminus\Omega^+$, 
a random trajectory starting at 
$x_0$ will have the form
\bel{X-}
X(t, x_0, \omega)~=~S^-_{T^-(t,x_0,\omega)} (x_0)\,,\eeq
where the time $T^-$ along the trajectory is a random variable with distribution 
\bel{T-} {\P\bigl\{ T^-(t, x_0, \omega)<s\bigr\}}~=~\P\left\{ s+
\sum_{y_j\in S^*\cap [ S^-_s(x_0),\, x_0]}  Y_j (\omega)~>~t\right\}.\eeq
The transition probabilities are thus given by 
\bel{TP3}\bega{l}
\ds
P_t\bigl( x_0,\, \,]x_0, +\infty[\bigr) ~=~0,\\[3mm]
\ds P_t\bigl( x_0,\, \,[S_s^-(x_0), x_0]\bigr)~=~\P\left\{ s+\sum_{y_j\in S^*\cap 
[ S^-_s(x_0), \, x_0]}  Y_j (\omega)~\geq ~t\right\}.\enda\eeq
\item[(iv)]
To complete the definition, it remains to define the transition probabilities in the case
$x_0 \in\Omega^+\cap\Omega^-\subseteq\Omega^*$.     In this case,
by construction we have $x_0= z_k$ for some $k$.
We then define
the random variable 
$X(t, x_0,\omega)$ by setting
\bel{X**} 
X(t, x_0,\omega)~=~\left\{ \bega{rl} S^+_{T^+(t, x_0, \omega)}(x_0)\quad
&\hbox{if}\quad Z_k(\omega)=1\,,\\[4mm]
S^-_{T^-(t, x_0, \omega)}(x_0)\quad
&\hbox{if}\quad Z_k(\omega)= -1\,.\enda\right.\eeq
By (\ref{pxj}), its distribution satisfies
\bel{TP4}P_t(x_0,A)~=~
\Theta(x_0)\cdot\P\left\{S^+_{T^+(t, x_0, \omega)}(x_0)\in A\right\} + (1-\Theta(x_0))\cdot 
\P\left\{S^-_{T^-(t, x_0, \omega)}(x_0)\in A\right\} \eeq
for every Borel set $A\subseteq\R$.
\endi

From the above construction it is clear that, in all cases (i)--(iv), every sample path
 $t\mapsto X(t, x_0,\omega)$ is a Carath\'eodory solution of the  Cauchy
 problem  (\ref{ode})-(\ref{ic}). Indeed, we are only adding a countable number of waiting times, when the random trajectories reach one of the points $y_j\in \S^*\subseteq 
 f^{-1}(0)$. 
 We conclude by stating the second main result proved in \cite{BMN}.

\begin{theorem}\label{t:23} Let $f$ be a function satisfying {\bf (A1)-(A2)}.
The following statements are equivalent.
\begi
\item[(I)] The random variables $X(t,x_0,\omega)$ yield
a Markov process whose sample
paths are solutions to the ODE (\ref{ode})-(\ref{ic}).
\item[(II)] There exist: (i) a positive, atomless Borel measure $\mu$ supported on $f^{-1}(0)$,
(ii) a countable set $\S\subseteq f^{-1}(0)$ of stationary points, (iii) a countable set $\S^*
= \{y_j: j\geq 1\}\subseteq f^{-1}(0)$ and corresponding numbers $\lambda_j>0$
determining the Poisson waiting times,
and (iv) a map $\Theta:\Omega^*\mapsto [0,1]$, such that
the  transition 
kernels $P_t(x_0,A)= \P\{ X(t,x_0,\omega)\in A\}$ coincide with the corresponding ones
constructed 
at (\ref{TP1}), (\ref{TP2}), (\ref{TP3}), (\ref{TP4}).
\endi \end{theorem}

\section{Approximating a deterministic semigroup by smooth flows}
\label{s:3}
\setcounter{equation}{0}

In this section we give a proof of
Theorem~\ref{t:1}, in several steps.
\v
{\bf 1.}  In analogy with Definition~\ref{d:increase}, we introduce:
\begin{definition}\label{d:incdyn} Given the semigroup $S:\R\times\R_+\mapsto\R$, we say that 
\begin{itemize}
\item an  interval $J= [a,b[$ or $J=\,]a,b[$ (open to the right) is a 
{\bf domain of increasing dynamics}
if, for every $x_0, \Hat x\in J$ with $x_0<\Hat x$ there exists $t>0$ such that 
$\Hat x ~=~S_t x_0$;
\item  an  interval $J= ]a,b]$ or $J=\,]a,b[$ (open to the left) is a 
{\bf domain of decreasing dynamics}
if, for every $x_0, \Hat x\in J$ with $\Hat x < x_0$ there exists $t>0$ such that 
$\Hat x ~=~S_t x_0$.
\end{itemize}
\end{definition}
We observe that, if two  domains of increasing (decreasing) dynamics $J, J'$ have nonempty
intersection, then the union $J\cup J'$ is also an interval of increasing (decreasing) dynamics.
We can thus partition the real line as
\bel{lp}
\R~=~\Omega_S^+\cup\Omega_S^-\cup\Omega_S^0\eeq
where
\begi
\item $\Omega_S^+=\ds\bigcup_{i\in \J^+}  J_i^+$ is the  union of countably many 
disjoint intervals $J_i^+= [a_i, b_i[$ or $J_i^+= \,]a_i, b_i[$ open to the right, where the dynamics is strictly increasing and  $b_i$ cannot be crossed from the left, i.e. $S_t(b_i-\ve)\leq b_i$ for all $t\geq 0$, $\ve\in ]0,b_i-a_i[$.
\item $\Omega_S^-=\ds\bigcup_{i\in \J^-}  J_i^-$ is the  union of countably many 
disjoint intervals $J_i^-=\,]c_i, d_i]$ or $J_i^-=\,]c_i, d_i[\,$ open to the left, where the dynamics is strictly decreasing and  $c_i$ cannot be crossed from the right, i.e. $S_t(c_i+\ve)\geq c_i$ for all $t\geq 0$, $\ve\in ]0,d_i-c_i[$.

\item $\Omega_S^0\subseteq f^{-1}(0)$ is the set of points $x_0$ such that 
$S_t(x_0)=x_0$ for all $t>0$.
\endi

Let $\mu$ be the positive atomless measure  supported on the 
set $f^{-1}(0)$, described at {\bf (Q1)}. Recall that $\mu$ is finite in any compact subset of $J^+_i$ or $J_i^-$.
\v
{\bf 2.} Let  $\ve>0$ be given.  As an intermediate step, we shall 
approximate $f$ with a piecewise constant function 
$g_\ve$.  
Consider the two  sets of indices 
\bel{Jep}
\J_\ve^+~\doteq~\bigl\{i\in \J^+\,;~ |J^+_i|=b_i-a_i \geq \ve\bigr\}, \qquad 
\qquad \J_{\ve}^-~\doteq~\bigl\{i\in\J^-\,;~ |J^-_i|= d_i-c_i\geq \ve\bigr\},
\eeq
labelling the intervals with length $\geq \ve$.
For each $i\in\J^+_{\ve}$, on the interval $J_i^+$ where the dynamics is increasing we perform the following construction.

CASE 1: If $J^+_i= ]a_i, b_i[$ is bounded and open, we insert points
$a_i< x_{i,1}<\dots< x_{i,N_i}<b_i $ so that
\bel{nod1}{\ve\over 3}<x_{i,1}-a_i<\ve,\qquad {\ve\over 3}<b_i-x_{i,N_i} <\ve,\qquad x_{i,k}-x_{i,k-1}~<~\ve\eeq
\bel{nod11}\tau_{i,k}~\doteq~\int_{x_{i,k-1}}^{x_{i,k}} {dy\over f(y)} + \mu\big(\big[x_{i,k-1}, x_{i,k}\big]\big)~<~\ve\quad\forall
k\in \{2,\ldots,N_i\}.\eeq
We then define the piecewise constant function $g_\ve: \Big]a_i+\ds{\ve\over 3}\,,~b_i-{\ve\over 3}\Big[\mapsto\R$ by setting
$$g_\ve(x)~=~0\qquad\forall x\in \Big]a_i+{\ve\over 3}, \, x_{i,1}\Big[\,
\cup\, \Big[x_{i,N_i}, \, b_i-{\ve\over 3}\Big[,$$
$$g_\ve(x)~=~{x_{i,k}-x_{i,k-1}\over \tau_{i,k}}
\qquad\hbox{if} \quad x\in\big[x_{i,k-1}, x_{i,k}\big[\quad\text{for some }k\in
\{2,\ldots,N_i\}\,.$$
CASE 2:  If  $J^+_i= [a_i, b_i[$ is bounded and half-open, we insert points
$a_i= x_{i,0}< x_{i,1}<\dots< x_{i,N_i}<b_i $ so that
\bel{nod2} {\ve\over 3}<b_i-x_{i,N_i} <\ve,\qquad x_{i,k}-x_{i,k-1}<\ve\eeq
\bel{nod21}\tau_{i,k}~\doteq~\int_{x_{i,k-1}}^{x_{i,k}} {dy\over f(y)} + \mu\big(\big[x_{i,k-1}, x_{i,k}\big]\big)<\ve\quad\forall
k\in\{1,\ldots,N_i\}\,.\eeq
We then define the piecewise constant function $g_\ve: \Big[\ds a_i-{\ve\over 3}\,,~b_i-{\ve\over 3}\Big[\mapsto\R$ by setting
$$g_\ve(x)~=~0\qquad\forall x\in \Big[x_{i,N_i}, b_i-{\ve\over 3}\Big[,$$
$$g_\ve(x)~=~{x_{i,k}-x_{i,k-1}\over \tau_{i,k}}
\qquad\hbox{if} \quad x\in\big[x_{i,k-1}, x_{i,k}\big[\quad\text{for some }k\in\{2,\ldots,N_i\}\,,$$
$$g_\ve(x)~=~{x_{i,1}-x_{i,0}\over \tau_{i,1}}
\qquad\hbox{if} \quad x\in\Big[a_i-{\ve\over 3},x_{i,1}\Big[\,.$$
If $J^+_i$ is unbounded,   we can perform  the same construction as above by inserting countably many points. In this case, $N_i=\infty$.
\v

Next, for each $i\in\J^-_{\ve}$, we perform an entirely similar construction on 
the  interval $J_i^-= ]c_i, d_i[$ 
 or $J_i^-= ]c_i, d_i]$ where the dynamics is decreasing. Notice that the choice (\ref{nod1}), (\ref{nod2}) 
 of the nodes $x_{i,k}$ of each partition imply that the 
 domains where the functions $g_\ve$ are defined are strictly disjoint. 
 
 We then extend the function $g_\ve$ to the whole real line $\R$, by defining
$g_\ve(x)~=~0$ on the remaining set.
 \v
 {\bf 3.} In this step, we approximate $g_\ve$
 with a smooth function $f_\ve :\R\mapsto\R$, with the following properties.
\begi
\item[(i)] Introducing the convex valued multifunction 
$G_\ve(x)\doteq \hbox{co}\{g_\ve(x+), \, g_\ve(x-)\}$, one has
$\mathrm{Graph}(f_{\ve})\subseteq B\left(\mathrm{Graph}(G_{\ve}),\ve\right)$.
Moreover,
\bel{nod113}f_\ve(x)~=~0\quad\forall
x\in\R\setminus\left[\Big(\bigcup_{i\in \J^+_{\ve}}\big[a_i-{\ve\over 2},\, 
b_i\big[\Big)\cup\Big(\bigcup_{i\in \J^-_{\ve}}\big]c_i,\,
d_i+{\ve\over 2}\big]\Big)\right].\eeq
\item[(ii)] For every $i\in\J_{\ve}^{+}$,  if $J_i^+=\,]a_i, b_i[$ is open then $f_{\ve}(a_i)=0$ and
\bel{nod111}
 \int_{x_{i,k-1}}^{x_{i,k}} {dy\over f_\ve^+(y)} ~=~\tau_{i,k}\quad\forall k=2,\ldots,N_i;
\eeq
On the other hand, if $J_i^+=[a_i, b_i[$ is half open, then (\ref{nod111}) holds for all $k=1,\dots,N_i$. 

\item[(iii)] The analogous conditions hold for the intervals $J_i^-$ with $i\in \J^-_{\ve}$.\endi
Notice that all the above properties can be achieved by a suitable mollification.
\v
{\bf 4.} We claim that, as $\ve\to 0+$, the functions  $f_{\ve}$ converge to $f$ in the sense of the graph. Since $\mathrm{Graph}(f_{\ve})\subseteq B\left(\mathrm{Graph}(G_{\ve}),\ve\right)$, it will be sufficient to show that
\bel{graph-ve}
(x,g_{\ve}(x))~\in~B\left(\mathrm{Graph}(F),{4\ve\over 3}\right)\qquad\forall x\in\R,
\eeq
since this yields $\mathrm{Graph}(G_{\ve})~\subseteq~B\left(\mathrm{Graph}(F),\ds{4\ve\over 3}\right)$. Let us prove (\ref{graph-ve}) for  $g_{\ve}(x)\neq 0$.
\begin{itemize}
\item If $x\in \Big[\ds a_i-{\ve\over 3},b_i\Big[$ for some $i\in\J^+_{\ve}$ then 
\[
x~\in~ \Big[\ds a_i-{\ve\over 3},x_{i,1}\Big[ \qquad\mathrm{or}\qquad x\in\big[x_{i,k-1}, x_{i,k}\big[\quad\text{for some }k\in\{2,\ldots,N_i\}.
\]
Assume that $x\in \Big[\ds a_i-{\ve\over 3},x_{i,1}\Big[$. Then 
\[
g_{\ve}(x)~=~{1\over \tau_{i,1}}\cdot \int_{0}^{\tau_{i,1}}f(S_{s}(x_{i,0}))ds~\in~\left[\inf_{y\in [x_{i,0},x_{i,1}]}f(y),\sup_{y\in [x_{i,0},x_{i,1}]}f(y)\right].
\]
Thus, since $0<x_{i,1}-a_i<\ve$, it follows  that $(x,g_{\ve}(x))\in B\left(\mathrm{Graph}(F),\ds{4\ve\over 3}\right)$. Otherwise,  if $x\in\big[x_{i,k-1}, x_{i,k}\big[$ for some $k\in\{2,\ldots, N_i\}$, then 
\[
g_{\ve}(x)~=~{1\over \tau_{i,k}}\cdot \int_{0}^{\tau_{i,k}}f(S_{s}(x_{i,k-1}))ds~\in~\left[\inf_{y\in \big[x_{i,k-1},x_{i,k}\big]}f(y),\sup_{y\in\big[x_{i,k-1},x_{i,k}\big]}f(y)\right],
\]
and this implies that $(x,g_{\ve}(x))\in \left(\mathrm{Graph}(F),\ve\right)$.
\item Similarly, one can show that (\ref{graph-ve}) holds  if $x\in \Big]c_i,d_i+\ds{\ve\over 3}\Big]$ for some $i\in\J^{-}_{\ve}$.
\end{itemize}

Let us now prove (\ref{graph-ve}) in the case where  $g_{\ve}(x)= 0$. 
This is trivial when $x\in\Omega_S^0$. 
On the other hand, if $x\notin \Omega_S^0\cup\left(\bigcup_{i\in\J^+_{\ve}}I_i^+\right)\cup\left(\bigcup_{i\in\J^-_{\ve}}I_i^-\right)$, then 
\[
|x- c_i|~<~\ve\quad\mathrm{for~some}~i\notin \J^+_{\ve}\qquad\mathrm{or}\qquad |x-b_i|~<~\ve\quad\mathrm{for~some}~i\notin \J^-_{\ve}.
\]
Recalling that $f(c_i)=f(b_i)=0$, we obtain (\ref{graph-ve}). Now, assuming that $x\in J^+_{i}$ for some $i\in\J^+_{\ve}$, we have  
\[
|x-b_i|~\leq~{4\over 3}\cdot\ve\qquad\mathrm{or}\qquad |x-a_i|~\leq~\ve,~~a_i\notin J_{i}^+.
\]
If $|x-b_i|\leq\ds {4\over 3}\ve$ then  (\ref{graph-ve}) holds because $f(b_i)=0$. 
Otherwise, since $a_i\notin J_{i}^+$, one has  $f(a_i)=0$ or $f(a_i-)<0<f(a_i+)$ and this also implies (\ref{graph-ve}). 

Similarly, one can show that   (\ref{graph-ve}) holds if $x\in J^-_i$ for some 
$i\in \mathcal{J}^-_{\ve}$. 
\medskip

{\bf 5.}  To complete the proof, it remains to show that, for each $x_0\in\R$, the unique solution $t\mapsto  S^{\ve}_t x_0$  to the 
Cauchy problem 
\[
\dot x ~=~ f_{\ve}(x),\qquad\qquad x(0)~=~ x_0,
\]
converges uniformly to $t\mapsto  S_t x_0$, i.e. 
\bel{unfc}
\lim_{\ve\to 0+}\left(\sup_{t\geq 0}\big|S^{\ve}_t(x_0)-S_t(x_0)\big|\right)~=~0.
\eeq
Several cases will be considered:
\v

CASE 1: $x_0\in J_i^+$, where $i\in \J^+_{\ve}$ and $J_i^+=\,]a_i, b_i[$ is open. 
For every $\ve$ such that  $$0~<~\ve~< ~\min \left\{x_0 - a_i, ~b_i-x_0,~{3\over 7}
\right\},$$ one then has
$$
\|f_{\ve}\|_{\infty}~\leq~\|f\|_{\infty}+{7\ve\over 3}~\leq~M+1\qquad\mathrm{and}\qquad x_0\in\big[x_{i,k-1}, x_{i,k}\big[~~\mathrm{for~some}~ k\in \{2,\dots, N_i\}.$$
Let $\tau_{0},\tau^{\ve}_0\in\,]0,\ve[$ be such that $S_{\tau_0}(x_0)=S^{\ve}_{\tau^{\ve}_0}(x_0)=x_{i,k}.$
 For any $0<t<\max\{\tau_{0},\tau^{\ve}_0\}$ one has
\[
\big|S^{\ve}_t(x_0)-S_{t}(x_0)\big|~\leq~(M+1)\ve.
\]
On the other hand, if $t\geq \max\{\tau_{0},\tau^{\ve}_0\}$, we have 
\begin{eqnarray*}
\big|S^{\ve}_t(x_0)-S_t(x_0)\big|&=&\big|S^{\ve}_{t-\tau^{\ve}_0}\big(x_{i,k}\big)-S_{t-\tau_0}\big(x_{i,k}\big)\big|\\
&\leq&\big|S^{\ve}_{t-\tau^{\ve}_0}\big(x_{i,k}\big)-S^\ve_{t-\tau_0}\big(x_{i,k}\big)\big|+\big|S^{\ve}_{t-\tau_0}\big(x_{i,k}\big)-S_{t-\tau_0}\big(x_{i,k}\big)\big|\\
&\leq&(M+1)\ve+\big|S^{\ve}_{t-\tau_0}\big(x_{i,k}\big)-S_{t-\tau_0}\big(x_{i,k}\big)\big|.
\end{eqnarray*}
To estimate the second term $\big|S^{\ve}_{t-\tau_0}\big(x_{i,k}\big)-S_{t-\tau_0}\big(x_{i,k}\big)\big|$, we consider two subcases: 
\begi
\item if $S_{t-\tau_0}\big(x_{i,k}\big)\in\big[x_{i,k_0-1}, x_{i,k_0}\big[$ for some $k_0\in
\{ k+1,\dots,N_i\}$, then 
\begin{eqnarray*}
\big|S^{\ve}_{t-\tau_0}\big(x_{i,k}\big)-S_{t-\tau_0}\big(x_{i,k}\big)\big|&=&\left|S^{\ve}_{t-\tau_0-\sum_{k'=k+1}^{k_0-1}\tau_{i,k'}}\big(x_{i,k_0-1}\big)-S_{t-\tau_0-\sum_{k'=k+1}^{k_0-1}\tau_{i,k'}}\big(x_{i,k_0-1}\big)\right|\\
&\leq&x_{i,k_0}-x_{i,k_0-1}~\leq~\ve\,; 
\end{eqnarray*}
\item Otherwise, if $S_{t-\tau_0}\big(x_{i,k}\big)\in\big[x_{i,N_i}, b_i\big]$ then 
\begin{eqnarray*}
\big|S^{\ve}_{t-\tau_0}\big(x_{i,k}\big)-S_{t-\tau_0}\big(x_{i,k}\big)\big|&=&\left|S^{\ve}_{t-\tau_0-\sum_{k'=k+1}^{N_i}\tau_{i,k'}}\big(x_{i,N_i}\big)-S_{t-\tau_0-\sum_{k'=k+1}^{N_i}\tau_{i,k'}}\big(x_{i,N_i}\big)\right|\\
&\leq&b_i-x_{i,N_i}~<~\ve.
\end{eqnarray*}
\endi
In  all cases we conclude
\bel{conv}
\big|S^{\ve}_t(x_0)-S_t(x_0)\big|~\leq~(M+2)\ve\quad \forall\,t\geq0\,.
\eeq
\v
CASE 2: $x_0\in J_i^+$, where $i\in \J^+_{\ve}$ and $J_i^+=[a_i, b_i[$ is half-open.    
Choosing $\ve<b_i-x_0$, we can obtain \eqref{conv} in a similar way as in Case 1.

CASE 3: $x_0\in J_i^-$, with $i\in \J^-_{\ve}$. The estimates are analogous to the previous ones.

CASE 4: $x_0\in \Omega_S^0\bigcup\Big(\bigcup_{i\notin \J^+_{\ve}}I_i^+\Big)\bigcup\Big(\bigcup_{i\notin \J^-_{\ve}}I_i^-\Big)$. In this case, we have that
\bel{sma}
|S_t(x_0)-x_0|~<~\ve~~~\qquad\forall t\geq 0.
\eeq
If $x_0=a_i$ for some $i\in \J^+_{\ve}$ or $x_0=d_i$ for some $i\in \J^-_{\ve}$, then $S^{\ve}_t(x_0)=x_0=S_t(x_0)$ for any $t,\ve>0$.
Otherwise, choosing $\ve>0$ such that 
$$x_0\notin \Big[ a_i -{\ve\over 2}\,,~a_i\Big]\quad\forall\,i\in \J^+_{\ve}\quad\text{and}\quad 
x_0\notin \Big[ d_i \,,~d_i+ {\ve\over 2}\Big]\quad\forall\,i\in \J^-_{\ve},$$
we obtain  $S^{\ve}_t(x_0)=x_0$ for all $t\geq 0$. Hence (\ref{sma}) yields (\ref{unfc}).
\endproof

\v
\section{Markov semigroups with a simpler dynamics}
\setcounter{equation}{0}
\label{s:4}

According to Theorem~\ref{t:23}, a Markov semigroup compatible with 
the ODE (\ref{ode}) 
is uniquely determined by the positive, atomless measure $\mu$ supported on $f^{-1}(0)$,
the countable sets $\S$, $\S^*$, and the maps $\Lambda:\S^*\mapsto \R_+$,
$\Theta:\Omega^*\mapsto [0,1]$, following the construction in Section~\ref{s:2}. 

The next lemma shows that this random dynamics can be approximated
with a simpler one, where the function $f$ is piecewise constant,
the measure $\mu$ vanishes, and the sets $\S,\S^*,\Omega^*$ contain finitely many points.
Here and throughout the following, on the set of Lipschitz functions $\vp:\R\mapsto\R$,
we use the norm
\bel{lipn}\|\vp\|_{Lip}~=~\|\vp\|_{W^{1,\infty}}~=~\max\left\{ \sup_x |\vp(x)|~,~\sup_{x\not= y}
{|\vp(x)-\vp(y)|\over |x-y|}\right\}.\eeq

\begin{lemma}\label{l:41} Let $f$ satisfy the assumptions {\bf (A1)-(A2)}, and let
$P_t= P_t(x_0,A)$ be the transition kernels for a Markov semigroup,
compatible with (\ref{ode}).    Then, for any given $T,R>0$, there exists a sequence of
transition kernels $P_t^n=P^n_t(x_0,A)$, such that
\bel{PPn} 
\lim_{n\to\infty}~\sup_{t\in [0,T], |x_0|\leq R}~\sup_{\|\vp\|_{Lip}\leq 1}
\left| \int \vp(x) P^n_t(x_0, dx)-\int \vp(x) P_t(x_0, dx)\right|~=~0.
\eeq
For each $n\geq 1$,  the transition kernels $P^n_t(\cdot,\cdot)$
define a Markov semigroup compatible with an ODE
\bel{oden}
\dot x~=~f_n(x),\eeq
where $f_n$ is piecewise constant with compact support.
The functions $f_n$ converge to $f$ in the sense of the graph.
Moreover, the corresponding sets $\S_n, \S^*_n, \Omega^*_n$ are finite and disjoint, 
while all the measures $\mu_n$ vanish.
\end{lemma} 


{\bf Proof.} For every fixed $N\ge 1$, consider the Markov process $X^N$
obtained by 
stopping the motion as soon as it exits from the interval $[-N, N]$. 
This is achieved by inserting a stopping time $\tau_N$:
$$X^N(t, x_0, \omega)~=~X\bigl(\tau_N, x_0,\omega\bigr),$$
where
$$\tau_N(x_0,\omega)~\doteq~\min\Big\{ t\geq 0\,;~~\bigl|X(t, x_0,\omega)\bigr|
\,\geq\,N\Big\}.$$
The corresponding  transition kernels $P^N_t (x_0,\cdot )$
define a Markov semigroup compatible with the ODE $\dot x = f^N(x)$, where
$$f^N(x)~=~\left\{ \bega{cl} f(x)\qquad &\hbox{if}~~|x|<N,\\[2mm]
0\qquad &\hbox{if}~~|x|\geq N.\enda\right.$$
For every given $x_0\in \R$ and  $t\geq 0$, as soon as $N\geq Mt+|x_0|$ we have the identity
\[
P^N_t(x_0,A)~=~P_t(x_0,A)\qquad\forall~\mathrm{Borel~set}~A\subseteq \R. \]
Here $M$ is the upper bound on $|f|$, introduced in the assumption {\bf (A1)}. Each $f_N$ has a compact support,  and moreover  the functions 
$f_N$ converge to $f$ in the sense of the graph. 

By replacing $f$ with $f^N$, it thus suffices to consider the case where
$f$ has compact support. The proof will be given in several steps.
\v
{\bf 1.} We begin by identifying the set of points where the dynamics stops forever:
$$\Omega^0_X~=~\bigl\{x\in \R\,;~~ \P\{X(1,x,\omega)=x\}=1\bigr\}~\subseteq~ f^{-1}(0).$$
We also introduce:
\begin{definition}\label{d:41}
We say that a
right-open interval $J^+=\,]a,b[$ or $J^+=[a,b[$ is an {\bf interval of increase} 
if
\[
\lim_{t\to\infty}\P\{X(t,x_1,\omega)>x_2\}~=~1\qquad\forall x_1,x_2\in J^+~~\mathrm{with}~~ x_1<x_2.
\]
Similarly, a left-open interval 
 $J^-=\,]c,d[$ or $J^-=]c,d]$ is an {\bf interval of decrease} 
if
\[
\lim_{t\to\infty}\P\{X(t,x_1,\omega)<x_2\}~=~1\qquad\forall x_1,x_2\in J^-_k~~\mathrm{with}~~ x_1>x_2.
\]
\end{definition}

If $J_1^+, J_2^+$ are two intervals of increase, with $J_1^+\cap J_2^+\not= \emptyset$, 
then the union $J_1^+\cup J_2^+$ is also an interval of increase. 
The same of course is true for intervals of decrease.
We can thus identify countably many maximal intervals of increase
$J^+_k=\,]a_k,b_k[$ or $[a_k,b_k[$,  $k\in \J^+$, 
and  countably many maximal intervals of decrease
 $J^-_k=]c_k,d_k[$ or $\,]c_k,d_k]$, $k\in \J^-$. 
 
As in  Theorem \ref{t:23}, the set of points where the dynamics stops for a random waiting time is 
\[
\S^*~=~\bigl\{x\in \R\,;~~0<\P\{X(1,x,\omega)=x\}<1\bigr\}.
\]
We also consider the countable set of stationary points
\[
\S~=~\Big(\{a_k,b_k\,;~~k\in\J^+\}\cup \{c_k,d_k\,;~~k\in\J^-\}\Big)\cap \Omega^0_{X}\,.
\]

{\bf 2.} Next, given $\ve>0$, among all the maximal intervals of increase or decrease, we select
a finite number, chosen as follows:
\begi
\item[(i)] All maximal intervals whose  length satisfies 
$b_k-a_k>\ve$, or $d_k-c_k>\ve$.
\item[(ii)] All maximal intervals to the left $]z_k^-,z_k[$ or to the right $]z_k,z_k^+[$ of a point $z_k\in \Omega^*$,
from where both a decreasing and an increasing trajectory can initiate, whose  length satisfies $z_k^+-z_k^->\ve$.
\endi
It is clear that the above rules select a finite set of intervals $J_k$, $k=1,\ldots,N$.
We now define a new dynamics, setting

\bel{epdy}f_\ve(x)~=~\begin{cases}0 &\forall x~\notin ~\Omega^*\cup J_1\cup\cdots\cup J_N\,,\\
f(x) &\mathrm{otherwise}.
\end{cases}
\eeq
By construction, we have
\[
X^{\ve}(t, x_0,\omega)~=~x_0\qquad 
\forall t\geq 0,~~x_0~\notin ~\Omega^*\cup J_1\cup\cdots\cup J_N
\]
and
\bel{df1}
\bigl|X^{\ve}(t, x_0,\omega)-X(t, x_0,\omega)\bigr|~\leq~\ve\qquad\forall t\geq 0, x_0\in\R.
\eeq
Hence, calling $P_t^{\ve}$ the corresponding
transition probability kernels, for every given $(t,x_0)\in [0,\infty[\times \R$ and  every 
continuous function $\varphi$, we have 
\bel{kc1}\bega{l}\ds
\Big| E[\varphi(X^{\ve}(t,x_0,\omega))]-E[\varphi(X(t,x_0,\omega))]\Big|~\leq~\int_{\W}|\varphi(X^{\ve}(t,x_0,\omega))-\varphi(X(t,x_0,\omega))|~d\,\P(\omega)\\[3mm]
\quad \ds\leq~\sup\Big\{ |\varphi(s_1)-\varphi(s_2)|\,;~~s_1,s_2\in [x_0- tM, ~x_0+ tM],
~|s_1-s_2|\leq \ve\Big\}.
\enda
\eeq
Since $\vp$ is uniformly continuous on bounded intervals, as $\ve\to 0$ the right hand side
of (\ref{kc1}) approaches zero.
This implies  the convergence
\bel{PPP} 
\lim_{\ve\to 0} \int \vp(x) P^{\ve}_t(x_0, dx)~=~\int \vp(x) P_t(x_0, dx).
\eeq
Thanks to the previous arguments, we can  now assume that the dynamics is stationary
outside a finite number of intervals of increase or decrease.
As shown in Fig.~\ref{f:ode13},  three cases must be considered:

\begin{figure}[ht]
\centerline{\hbox{\includegraphics[width=12cm]{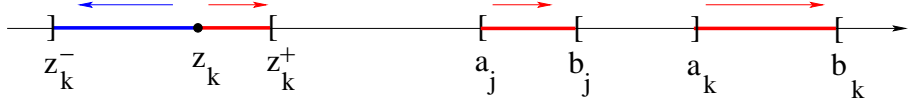}}}
\caption{\small  The three cases considered in the proof of Lemma~\ref{l:41}.
}
\label{f:ode13}
\end{figure}

Case 1: an open  maximal interval of increase $J_k\doteq\,]a_k, b_k[\,$, of length $b_k-a_k>\ve$.
In this case, we consider the smaller interval 
$$J_{k,\ve}~\doteq~\begin{cases} ]a_k, \, b_k-\ve/2[ &\mathrm{if}\qquad \mu(]a_k,a_k+\ve/2[)~<~\infty\cr\cr
]a_k+\delta_{\ve}, \, b_k-\ve/2[ &\mathrm{if}\qquad \mu(]a_k,a_k+\ve/2[)~=~ \infty
\end{cases}$$
where $0<\delta_{\ve}<\ve/2$ sufficiently small such that $\mu([a_k+\delta_{\ve},a_{k}+\ve])> T$ and so $X(T,x_0,\omega)<a_{k}+\ve$ for every $x_0\in ]a_k,a_k+\delta_{\ve}[$.

Case 2: a half-open maximal interval of increase $J_j\doteq [a_j, b_j[\,$, of length $b_j-a_j>\ve$.
In this case, we consider the smaller interval $J_{j,\ve}\doteq [a_j, \, b_j-\ve/2[\,$.

Case 3: a point $z_k\in \Omega^*$, from where both increasing and decreasing trajectories can originate, with probability $\theta\in [0,1]$ and $1-\theta$, respectively.
We then call $]z_k^-, z_k[$ and $]z_k, z_k^+[$ the corresponding maximal intervals
of decrease and of increase, to the left and to the right of $z_k$ such that $z_k^+-z_k^->\ve$.
In this case, we set $J_k\doteq\, ]z_k^-, z_k^+[\,$ and consider the smaller interval
$J_{k,\ve}\doteq \,]z_k^- + \ve(z_k-z_k^-), ~z_k^+ - \ve(z_k^+-z_k)[\,$.
\v
{\bf 3.} In all three cases, we modify the dynamics by setting
\bel{stok}
f_\ve(x)~=~0\qquad \hbox{if}\quad x\in J_k\setminus J_{k,\ve}\,.\eeq
In other words, the dynamics is stopped outside $J_{k,\ve}$.
This implies
\bel{edy1}
X^\ve(t,x_0, \omega)~=~x_0\quad\hbox{for all} ~~x_0\in\, J_k\setminus J_{k,\ve}\,.\eeq
Moreover, for $x_0\in J_{k,\ve}$, in Cases 1 and 2 we have
$$X^\ve(t,x_0, \omega)~=~\min\bigl\{ X(t, x_0,\omega), \,b_k-\ve/2\bigr\}\,,$$
while in Case 3:
$$X^\ve(t,x_0, \omega)~=~\left\{
\bega{rl}   X(t, x_0,\omega)\qquad &\hbox{if} \quad X(t, x_0,\omega)\in 
[z_k^- + \ve(z_k-z_k^-)\,,~z_k^+ -\ve(z_k^+-z_k)]\,,\\[3mm]
z_k^- + \ve(z_k-z_k^-)&\hbox{if} \quad X(t, x_0,\omega)~\leq~\,z_k^- + \ve(z_k-z_k^-),\\[3mm]
z_k^+ -\ve(z_k^+-z_k)&\hbox{if} \quad X(t, x_0,\omega)~\geq~\,z_k^+ -\ve(z_k^+-z_k).\enda
\right.$$
Since the new dynamics is stationary on each $J_k\setminus J_{k,\ve}$, 
the set $\S^*$ of points where a random waiting time occurs, as well as the measure $\mu$, 
can now be replaced by $\S^*_\ve$ and $\mu_\ve$ respectively.
Here
\bel{Smue}\S^*_\ve~\doteq~\S^*\cap \Big(\bigcup_k J_{k,\ve}\Big),\qquad
\qquad \mu_\ve(A)~\doteq~\mu \left( A\cap \Big(\bigcup_k J_{k,\ve}\Big)\right),\eeq
for every Borel set $A$.

We observe that, for a fixed $x_0$, we have
\bel{df3}
\bigl|X^{\ve}(t, x_0,\omega)-X(t, x_0,\omega)\bigr|~\leq~\ve\qquad\forall t\in [0,T], \omega\in\mathcal{W}
\eeq
Therefore, the same argument used at (\ref{kc1}) implies the convergence
(\ref{PPP}).
\v
{\bf 4.} For a fixed $\ve>0$, consider the approximate dynamics constructed in the 
previous step.    We claim that the measure $\mu^\ve$  in (\ref{Smue}) 
has finite total mass.

Indeed, consider an interval $J_k= \,]a_k, b_k[\,$ with $\mu(]a_k,a_k+\ve/2[)=\infty$, as in a subcase of Case 1.
By assumption there are some random trajectories 
starting at $a_k+\delta_{\ve}$ and reaching $b_k-\ve/2$ in finite time.
But, for every trajectory, the time needed to cross the whole interval $J_{k,\ve}$
is $\geq \mu(J_{k,\ve})$.    Hence this value must be finite.
The other two cases are entirely similar.

Since the intervals $J_{k,\ve}$ are finitely many, this implies 
$$\mu^\ve(\R)~=~
\sum_{k=1}^N\mu(J_{k,\ve})~<~+\infty.$$

{\bf 5.} Next, by a further approximation, we can modify the dynamics so 
that  each interval $J_{k,\ve}$, $k=1,\ldots,N$, contains at most finitely many points 
$y_j\in \S^*$ and thus the set $\S^*$ is finite. 
To justify this, let us consider a maximal  interval of increase $J_{k,\ve}$, as in Case 1.
By construction, there is a positive probability of moving from $a_k+\ve/2$ to $b_k-\ve/2$
in finite time. Hence, for some $t>0$ one has
\[
\eta~\doteq~\P\left\{X(t,a_k+\ve/2,\omega)\leq b_k-\ve/2\right\}~<~1\,.
\]
By the Markov property, we have 
\[
\P\left\{X(mt,a_k+\ve/2,\omega)<b_k-\ve/2\right\}~\leq~\eta^m.
\]
This implies that the sum of all random waiting times $Y_j\in J_{k,\ve}$ satisfies
\[
\P\left\{\sum_{y_j\in [a_k+\ve/2,b_k-\ve/2]}Y_j\geq mt\right\}~\leq~\eta^m\qquad\forall m\geq 1.
\]
In particular, 
$$E\left[\sum_{y_j\in [a_k+\ve/2, b_k-\ve/2[} Y_j\right]~=~\sum_{y_j\in [a_k+\ve/2, b_k-\ve/2[} {1\over\lambda_j}~<~+\infty.$$
Given an auxiliary constant $\delta>0$, we  can now 
\begi
\item remove all waiting times $Y_j$ with 
$$y_j\in [a_k, a_k+\ve/2]\cup [b_k-\ve/2, b_k],$$
\item remove all waiting times $Y_j$, with $y_j\in J_{k,\ve/2}= \, ]a_k+\ve/2,\, b_k-\ve/2[\,$ and 
$j> N_k$, choosing $N_k$ large enough so that
\bel{Edel}E\left[\sum_{y_j\in J_{k,\ve}, ~j> N_k} Y_j\right]~=~\sum_{y_j\in   J_{k,\ve}, ~j>N_k} {1\over\lambda_j}~<~\delta.\eeq
\endi
This will produce a new family of random trajectories $\Tilde X^\ve(\cdot, x_0, \omega)$.
Since we are removing some of the waiting times, within the interval $J_{k,\ve}$
the new trajectories will be shifted forward,  by an amount 
\bel{TXe}0~\leq~\Tilde X^\ve(t, x_0, \omega)-X^\ve(t, x_0, \omega)~\leq ~M\!\cdot\! \left(\sum_{y_j\in   J_{k,\ve}, ~j>N_k} Y_j(\omega)\right),\eeq
where $M$ is the upper bound on $f$, introduced in {\bf (A1)}.

For any given $\delta>0$, we now set 
$$
\W_\delta~\doteq~\left\{\omega\in\W\,;~~M\!\cdot\! \left(\sum_{y_j\in   J_{k,\ve}, ~j>N_k} Y_j(\omega)\right)\geq \sqrt{\delta}\right\}.
$$
By Chebyshev's inequality, it follows
\[
\P(\W_{\delta})~\leq~{M\over \sqrt{\delta}}\cdot E\left[ \left(\sum_{y_j\in   J_{k,\ve}, ~j>N_k} Y_j(\omega)\right)\right]~\leq~M\sqrt{\delta}.
\]
Therefore, given any $x_0\in J_{k,\ve}$ and $t>0$,  and any bounded 
continuous function $\varphi$, 
using (\ref{TXe}) we estimate 
\bel{ETE}\bega{l}\ds
\Big| E[\varphi(\Tilde X^{\ve}(t,x_0,\omega))]-E[\varphi(X(t,x_0,\omega))]\Big|
~\leq~\int_{\W}
\left|\varphi(\Tilde X^{\ve}(t,x_0,\omega))-\varphi(X^\ve(t,x_0,\omega))\right|
~d\P(\omega)\\[4mm]
\qquad \ds\leq~2\|\varphi\|_{\infty}\cdot \P(\W_\delta)+\int_{\W\backslash\W_\delta } 
\left|\varphi(\Tilde X^{\ve}(t,x_0,\omega))-\varphi(X^\ve(t,x_0,\omega))\right|
~d\P(\omega)\\[4mm]
\qquad \ds 
\leq~2M\|\varphi\|_{\infty}\!\cdot \!\sqrt{\delta}+\sup_{s_1,s_2\in J_{k,\ve},~|s_1-s_2|
\leq \sqrt{\delta}}~
\bigl|\varphi(s_1)-\varphi(s_2)\bigr|.
\enda \eeq
By choosing $\delta>0$ in (\ref{Edel}) sufficiently small, 
the difference of the expectations in (\ref{ETE}) can be rendered as small as we like.

The same ideas apply to intervals  such as $J_k= [a_k, b_k[\,$   or $\,]z_k^-, z_k^+[\,$, 
described in Cases 2 and 3 respectively. Of course, intervals of decrease can be 
handled in the same way.  
\v

{\bf 6.} In general, the semigroup may contain an exponential waiting time
also at points $z_k\in \Omega^*$ where trajectories can move both 
upward or downward, with given probabilities. 
However, we can construct an approximation so that the two sets
$\S^*$ and $\Omega^*$ are disjoint. 

Indeed, consider a point in the intersection: $z_k=y_j\in \Omega^*\cap \S^*$, 
for some indices $k,j$. Since $\S^*$ is finite, we can choose $\ve>0$ small enough
so that the interval 
$$[y'_j, y''_j]~\doteq~ [y_j-\ve\,,~ y_j+\ve]$$
does not contain any other point of $S^*$.  We now redefine $f$ at these two points,
by setting $f(y'_j)=f(y''_j)=0$. 
Moreover,  we construct a new  family  of transition kernels $P^\ve_t(\cdot,\cdot)$ by removing the random 
waiting time $Y_j$ at the point  $y_j=z_k$, and inserting  two 
waiting times  $Y'_j=Y''_j=Y_j$,
at the two points $y'_j$ and $y''_j$. 
In this case, for each random trajectory we have
\[
\bigl|X^{\ve}(t, z_k,\omega)-X(t, z_k,\omega)\bigr|~\leq~\max 
\bigl\{|y'_j-y_j|, |y''_j-y_j|\bigr\}~<~\ve
\]
and 
\[
X^{\ve}(t, x,\omega) ~=~ X(t, x,\omega)\qquad\forall x\notin [y_j-\ve,\, y_j+\ve],
\]
for every $t>0$, $\omega\in \W$. Therefore, as $\ve\to 0$,
the same argument used at (\ref{kc1}) implies the convergence in distribution.

Summarizing the previous analysis, we can 
approximate the original Markov semigroup with a new semigroup where:
\begi
\item[(i)] The dynamics is stationary outside finitely many intervals of increase or decrease.
\item[(ii)] The measure $\mu$ is finite.
\item[(iii)] The sets $\S, \S^*, \Omega^*$ are finite and  disjoint.
\endi

{\bf 7.} To achieve the proof, we still need to  approximate the dynamics with another Markov semigroup where $f$ is piecewise constant, and $\mu $ vanishes. 
To fix ideas, we will show how to modify  $f$ on  a maximal interval of increase $J_j$. The construction is entirely similar in the case of a maximal interval of decrease. 

Given $\ve>0$ sufficiently small, two cases will be considered: 

{\bf CASE 1:} The interval  $J_j = ]a_j,b_j[$ is open. 

By the same technique  used in Section~\ref{s:3}, for the proof of Theorem~\ref{t:1}, (step {\bf 2}, CASE 1),
 we can approximate $f$ by a function $f_\ve$ such that  
\begin{itemize}
\item[(i)] $f_{\ve}$ is piecewise constant and 
\[
\mathrm{Graph}(f_{\ve})~\subseteq~B(\mathrm{Graph}(F),\ve),\quad\qquad f_{\ve}(x)~=~f(x)\qquad\forall x\in \left(\R\backslash J_j\right)\cup \S^*;
\]
\item[(ii)] For every $[x_1,x_2]\subseteq  \bigl ]a_j,b_{j}-\ve/2\bigr]$, 
one has
\bel{d-time}
\left|\int_{x_1}^{x_2}{dy\over f(y)}+\mu([x_1,x_2])-\int_{x_1}^{x_2}{dy\over f_{\ve}(y)}\right|~\leq~\ve.
\eeq
\end{itemize}
We then take $\mu_\ve$ to be the  zero measure.

In view of (i)--(ii), for every  $x_0\in  \R\backslash J_j$ one has 
\bel{XXe}
X(t,x_0,\omega)~=~X^{\ve}(t,x_0,\omega)\qquad\forall (t,\omega)\in [0,\infty[\times \W.
\eeq
On the other hand, for a fixed $x_0\in J_j$,  
if $\min\{X(t,x_0,\omega),X^{\ve}(t,x_0,\omega)\}\geq b_j-\ve$,  we have 
\[
\bigl|X^{\ve}(t,x_0,\omega)-X(t,x_0,\omega)\bigr|~\leq~|(b_j-\ve)-
b_j|~=~\ve.
\]
Otherwise, without loss of generality,  assume that 
$$X(t,x_0,\omega)~>~X^{\ve}(t,x_0,\omega)\qquad\mathrm{and}\qquad X^{\ve}(t,x_0,\omega)~<~b_j-\ve.$$
Let $\tau^{\ve}(\omega), \tau_1(\omega),\tau(\omega)$ be such that 
\[
S^+_{\tau(\omega)}(x_0)~=~X(t,x_0,\omega)\qquad\mathrm{and}\qquad S^+_{\tau_1(\omega)}(x_0)~=~S^{+,{\ve}}_{\tau^{\ve}(\omega)}(x_0)~=~X^{\ve}(t,x_0,\omega)
\]
where $S^+$ and $S^{+,\ve}$ are defined  as in (\ref{S+x}). From (\ref{d-time}), we have 
\[
\left|\tau_1(\omega)-\tau^{\ve}(\omega)\right|~=~\left|\int_{x_0}^{S^+_{\tau_1(\omega)}(x_0)}{dy\over f(y)}+\mu([x_0,S^+_{\tau_1(\omega)}(x_0)])-\int_{x_0}^{S^+_{\tau_1(\omega)}(x_0)}{dy\over f_{\ve}(y)}\right|~\leq~\ve.
\]
Thus,
\begin{eqnarray*}
0~<~\tau(\omega)-\tau_1(\omega)&\leq&\left(t- \sum_{y_j\in \left[x_0,S^+_{\tau(\omega)}(x_0)\right[} Y_{j}(\omega)\right) - \tau^{\ve}(\omega)+\ve\\
&\leq&\left(t- \sum_{y_j\in \left[x_0,S^{+,\ve}_{\tau^{\ve}(\omega)}(x_0)\right]} Y_{j}(\omega)- \tau^{\ve}(\omega)\right)+\ve~\leq~\ve
\end{eqnarray*}
and this implies
\bel{et11}
|X^{\ve}(t,x_0,\omega)-X(t,x_0,\omega)|~=~\Big|S^+_{\tau(\omega)}(x_0)-S^+_{\tau_1(\omega)}(x_0)\Big|~\leq~M\ve.
\eeq

{\bf CASE 2:}  The interval  $J_j = [a_j,b_j[$ is half open.   

By the same technique use in Section~\ref{s:3}  (step {\bf 2}, CASE 2),
we can approximate $f$ by  a function $f_\ve$ such that  
\begin{itemize}
\item[(i)] $f_{\ve}$ is piecewise constant, and 
\[
\mathrm{Graph}(f_{\ve})~\subseteq~B(\mathrm{Graph}(F),\ve),\quad\qquad f_{\ve}(x)~=~f(x)\qquad\forall x\in \left(\R\backslash J_j\right)\cup \S^*;
\]
\item[(ii)] For every $[x_1,x_2]\subseteq \left [a_j,b_{j}-\ve/2\right]$, it holds
\bel{d-time2}
\left|\int_{x_1}^{x_2}{dy\over f(y)}+\mu([x_1,x_2])-\int_{x_1}^{x_2}{dy\over f_{\ve}(y)}\right|~\leq~\ve.
\eeq
\end{itemize}
Again, we take $\mu_\ve$ to be the zero measure.  

For every  $x_0\in  \R\backslash J_j$, the identity (\ref{XXe}) again holds.  On the other hand, if $x_0\in J_j$,   the same argument  used in the previous case yields (\ref{et11}).


In both cases, given any $x_0\in \R$, for each random trajectory we have
\[
\bigl|X^{\ve}(t, x_0,\omega)-X(t, x_0,\omega)\bigr|~\leq~M\ve\qquad\forall t\in [0,T],\omega\in \W.
\]
Therefore,  
the same argument used at (\ref{kc1}) implies the convergence (\ref{PPP}).
\endproof

\v
\section{An approximation lemma}
\setcounter{equation}{0}
\label{s:5}
The main goal of the next two sections is to prove Theorem~\ref{t:2}, 
showing that every Markov semigroup
whose trajectories are solutions to the ODE (\ref{ode}) can be approximated by a sequence of diffusion processes
with smooth coefficients (\ref{diffu}).  

Thanks to Lemma~\ref{l:41}, we can assume that the function $f$
is piecewise constant, the sets $\S, \S^*, \Omega^*$ are finite, 
and the measure $\mu$ vanishes.  
To fix the ideas, let 
\bel{x0N}x_0~<~x_1~<~x_2~<~\cdots~< x_N\eeq
be a list of all points in $\S, \S^*, \Omega^*$, together with all points where $f$ has a jump.
Consider the midpoints
\bel{y1N} y_1~<~y_2~<~\cdots~< y_N\,,\qquad\qquad y_j = {x_{j-1} + x_j\over 2}\,.\eeq
For any random trajectory starting at a point $\bar x$, say $s\mapsto X(s,\bar x, \omega)$ we define the stopping time
\bel{tom}\tau(\bar x, \omega)~=~\min\Big\{ s\geq 0\,;~X(s,\bar x,\omega)= y_j\quad\hbox{for some}~ y_j\not= \bar x\Big\}.\eeq
In other words, 
if $y_{j-1}<\bar x<y_j$, then $\tau$ is the first time when the random trajectory
hits either $y_{j-1} $ or $y_j$.    In the special case where $\bar x = y_k$, then 
$\tau$ is the first time when the trajectory hits either $y_{k-1}$ or $y_{k+1}$.

Similarly, given a sequence of Markov processes $X_n$, $n\geq 1$, we define 
\bel{tomn}\tau_n(\bar x, \omega)~=~\min\Big\{ s\geq 0\,;~X_n(s,\bar x,\omega)= y_j\quad\hbox{for some}~ y_j\not= \bar x\Big\}.\eeq

\begin{figure}[ht]
\centerline{\hbox{\includegraphics[width=8cm]{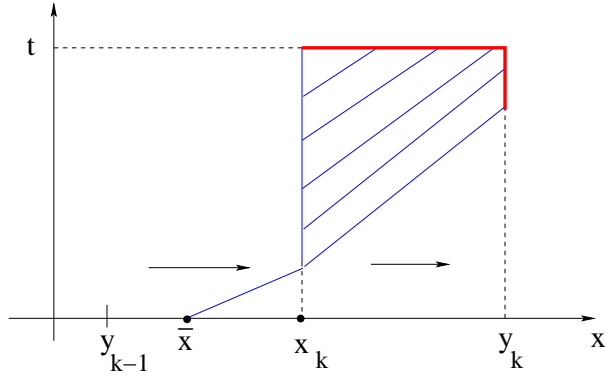}}}
\caption{\small Random trajectories of the Markov semigroup, starting at $\bar x$.
Here the speed $f(x)$ is constant for $x<x_k$ and $x> x_k$.    At $x_k$ there is a random waiting time.    The thick red lines denote the support of the measure $\mu^{t,\bar x}$.
}
\label{f:ode15}
\end{figure}


Throughout the following, we adopt the notation
$a\wedge b\doteq \min\{a,b\}$.
Using the above stopping times, given an initial point $\bar x$,
 for every  $t>0$ we consider the probability measures 
 $\mu^{t, \bar x}$, $\mu_n^{t, \bar x}$ on $\R_+\times \R$
defined as the push-forward of the maps
$$
{\omega~\mapsto~\bigl(\tau(\bar x, \omega)\wedge t, ~X(\tau(\bar x, \omega)\wedge t, \,\bar x, \omega)\bigr)},\qquad \quad{\omega~\mapsto~\bigl(\tau_n(\bar x, \omega)\wedge t, ~X_n(\tau_n(\bar x, \omega)\wedge t, \,\bar x, \omega)\bigr).}
$$
More precisely, for every open set $V\subset \R_+\times \R$, we define
\bel{muto}
{\mu^{t, \bar x}(V)~=~\meas\Big\{ \omega\,;~~\bigl(\tau(\bar x, \omega)\wedge t, ~X(\tau(\bar x, \omega)\wedge t, \bar x, \omega)\bigr)\in V\Big\}},
\eeq
\bel{muto-n}
{\mu_n^{t, \bar x}(V)~=~\meas\Big\{ \omega\,;~~\bigl(\tau_n(\bar x, \omega)\wedge t, ~X_n(\tau_n(\bar x, \omega)\wedge t, \bar x, \omega)\bigr)\in V\Big\}.}
\eeq
As shown in Fig.~\ref{f:ode15},
all measures $\mu^{t, \bar x}$ and $\mu_n^{t, \bar x}$ are supported within the set
\bel{spmu}\Sigma~\doteq~\{t\}\times \R ~\cup~ [0,t]\times \{y_1,y_2,\ldots, y_N\}.\eeq
The next lemma establishes a ``local to global" result.
If the random dynamics converge separately on each subinterval 
$[y_{j-1}, y_j]$ (i.e., after introducing the stopping times at each point $y_j$), 
then the weak convergence of the transition kernels (\ref{EP}) holds as well,
for every initial point $\bar x$ and every $t>0$.
\begin{lemma}\label{St2} Consider a Markov semigroup compatible with the ODE
(\ref{ode}), where $f$ is piecewise constant, the sets $\S, \S^*, \Omega^*$ are finite, 
and the measure $\mu$ vanishes.  
In addition, consider a sequence of  diffusions processes 
$X^{(n)}= X^{(n)}(s, \bar x,\omega)$ as in (\ref{diffu}), with 
$f_n\to f$ in the sense of the graph and $\sigma_n\downarrow 0$.
Assume that, for every $ \bar x\in \R$ and $t>0$, we have the weak convergence of the corresponding 
mesures:
\bel{wcon}
\mu^{t, \bar x}_n\wto \mu^{t, \bar x},\eeq
defined as in (\ref{muto})-(\ref{muto-n}).
Then the weak convergence  of the transition kernels (\ref{EP}) holds.
\end{lemma}

{\bf Proof.}  {\bf 1.} Given $\bar{x}\in \R$, for every $p=1,2,3,\ldots$, define the $p$-th random stopping time
by induction, in the obvious way.  
Namely, $\tau^1(\bar{x},\omega)$ is the first time when the trajectory $X(s,\bar{x},\omega)$ starting from $\bar{x}$ hits one of the points $y_j$.   
By induction, the $p$-th stopping time is
$$\tau^p(\bar{x},\omega)~=~\min\Big\{ s> \tau^{p-1}(\bar{x},\omega)\,;~~X(s,\bar{x},\omega) = y_j
\quad\hbox{for some }~y_j \not= X(\tau^{p-1}(\bar{x},\omega),\bar{x},\omega)\Big\}.$$
The corresponding measures $\mu^{t,\bar x}_p$ are defined as in (\ref{muto}),
replacing $\tau$ with $\tau^p$, i.e.,
\[
\mu^{t,\bar x}_p(V)~=~\meas\Big\{ \omega\,;~~\bigl(\tau^p(\bar x, \omega)\wedge t, ~X(\tau^p(\bar x, \omega)\wedge t, \bar x, \omega)\bigr)\in V\Big\}.
\]
The same construction can of course be repeated for each of the Markov processes $X_n$,
thus defining a sequence of measures $\mu^{t,\bar x}_{p,n}$, $n\geq 1$.
\medskip

{\bf 2.} We observe that, for a sequence of probability measures $\mu_n$ on the real line,
the weak convergence $\mu_n\wto\mu$ is equivalent to the convergence
of the distribution functions $F_n\to F$ in $\L^1_{loc}(\R)$, where
$$F_n(x)~=~\mu_n\bigl(\,]-\infty, x]\bigr),\qquad F(x)~=~\mu\bigl(\,]-\infty, x]\bigr)
\qquad\forall x\in \R.$$
In the case we are presently considering, 
since all of measures $\mu_{p,n}^{t,\bar x}$ 
are supported on the 1-dimensional set $\Sigma$ at (\ref{spmu}), the weak 
convergence
\bel{wco}\mu^{t,\bar x}_{p,n}~\wto~\mu^{t,\bar x}_{p}\eeq
is equivalent to the $\L^1$ convergence of a finite family of distribution functions.

More precisely, for every $j=1,\ldots, N$, consider the nondecreasing functions
\bel{Fdef} F^{t,\bar x}_{p,n,j} (s)~\doteq~\mu^{t,\bar x}_{p,n}\Big( [0,s]\times \{y_j\}\Big),
\qquad F^{t,\bar x}_{p,j}  (s)~\doteq~\mu^{t,\bar x}_{p}\Big( [0,s]\times \{y_j\}\Big),
\qquad s\in [0,t].\eeq
Moreover, for $y\in \R$, consider the nondecreasing functions
\bel{Gdef}\bega{l}\ds G^{t,\bar x}_{p,n}   (y)~\doteq~\mu^{t,\bar x}_{p,n}\Big( \{t\}\times ]-\infty, y]\Big)+
\sum_{y_j\leq y} \mu^{t,\bar x}_{p,n}\Big( [0,t]\times \{y_j\}\Big),\\[6mm]
\ds  G^{t,\bar x}_{p}   (y)~\doteq~\mu^{t,\bar x}_{p}\Big( \{t\}\times ]-\infty, y]\Big)+
\sum_{y_j\leq y} \mu^{t,\bar x}_{p}\Big( [0,t]\times \{y_j\}\Big)
.\enda\eeq
Then the weak convergence (\ref{wco}) is then equivalent to the 
convergence 
\bel{l1c}F^{t,\bar x}_{p,n,j} \to F^{t,\bar x}_{p,j} ~~\hbox{in}~~\L^1([0,t])
\qquad\hbox{and}\qquad G^{t,\bar x}_{p,n}\to G^{t,\bar x}_p~~\hbox{in}~~\L^1_{loc} (\R).
\eeq

\begin{figure}[ht]
\centerline{\hbox{\includegraphics[width=12cm]{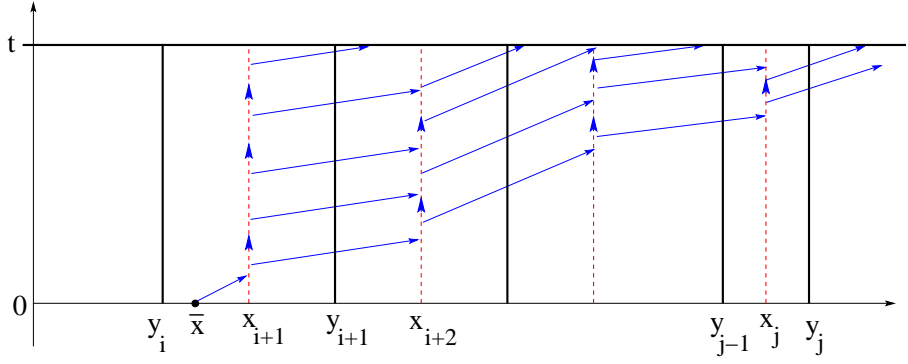}}}
\caption{\small The thick lines describe the set $\Sigma$ in (\ref{spmu}).  In this picture, 
trajectories
of the Markov semigroup start at $\bar x$ and move with positive speed, with 
random waiting times at the points $x_{i+1}, x_{i+2},\ldots$~    If $j = i+p$, the $p$-th stopping time
occurs when a trajectory reaches the point $y_j$, with $j= i+p$.
}
\label{f:ode16}
\end{figure}
\v
{\bf 3.} 
By induction, assuming the weak convergence 
\bel{wcm}\mu^{t,\bar x}_{p-1,n}~\wto~\mu^{t,\bar x}_{p-1}\,,\eeq
we need to prove the weak convergence (\ref{wco}).

To fix ideas, assume that all trajectories of the Markov semigroup starting from 
$\bar x$ are non-decreasing. Moreover, assume $y_i\leq \bar x < y_{i+1}$.
Setting $j = i+p$, this implies
\bel{spi}
\bega{rl}\ds 
Supp\bigl(\mu^{t,\bar x}_{p-1}\bigr)~\subseteq~\Sigma_{j-1}&\doteq~\{t\}\times \R ~\cup~ [0,t]\times \{y_{j-1}\}\,,\\[3mm]
Supp\bigl(\mu^{t,\bar x}_{p}\bigr)~\subseteq~\Sigma_{j}&\doteq~\{t\}\times \R ~\cup~ [0,t]\times \{y_{j}\}\,.\enda \eeq

We start by proving the convergence $F^{t,\bar x}_{p,n,j}\to F^{t,\bar x}_{p,j}$ in $\L^1\bigl([0,t])$.
Let $\ve>0$ be given.   By the inductive assumption, we have the convergence 
$F^{t,\bar x}_{p-1, n,j-1}\to F^{t,\bar x}_{p-1, j-1}$ in $\L^1\bigl([0,t]\bigr)$.

Call 
\bel{Gan}\bega{rl}\Gamma_n(\tau)&\doteq~\mu_{1,n}^{t, y_{j-1}}\bigl( [0,\tau]\times \{y_j\}\bigr),\\[3mm]
\Gamma(\tau)&\doteq~\mu_{1}^{t, y_{j-1}}\bigl( [0,\tau]\times \{y_j\}\bigr).
\enda\eeq
the probability that a random trajectory starting at $y_{j-1}$ reaches $y_j$ within time $\tau$
(without ever touching $y_{j-2}$, in the case of a diffusion).
The function $F^{t,\bar x}_{p,n,j}$
can now be computed by the Stieltjes integral 
\bel{Fj}
F^{t,\bar x}_{p,n,j}(\tau)~=~\int_0^\tau \Gamma_n(\tau-\zeta) \, d F^{t,\bar x}_{p-1,n,j-1}(\zeta).\eeq
We observe that all functions $\Gamma_n,\Gamma,  F^{t,\bar x}_{p-1,n,j-1}, F^{t,\bar x}_{p-1,j-1}$ are non-decreasing
with values in $[0,1]$. Moreover
$$\Gamma_n(0)~=~\Gamma(0)~=~ F^{t,\bar x}_{p-1,n,j-1}(0)~=~ F^{t,\bar x}_{p-1,j-1}(0)~=~0.$$
By the convergence
$\Gamma_n\to \Gamma$  and  $F^{t,\bar x}_{p-1,n,j-1}\to F^{t,\bar x}_{p-1,j-1}$  in $\L^1\bigl([0,t]\bigr)$
it thus follows the convergence of the integral functions in (\ref{Fj}), namely
$F^{t,\bar x}_{p,n,j}\to F^{t,\bar x}_{p,j}$  in $\L^1\bigl([0,t]\bigr)$.
\v

\begin{figure}[ht]
\centerline{\hbox{\includegraphics[width=8cm]{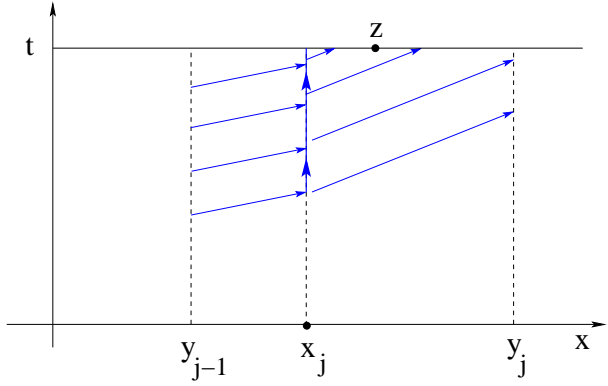}}}
\caption{\small Starting with the distribution $\mu^{t,\bar x}_{n,p-1}$ restricted to the
vertical segment $[0,t]\times \{y_{j-1}\}$, by the formula
(\ref{Lam}) one can compute the probability that trajectories of the semigroups reach a point 
$y\leq z$ at time $t$.
}
\label{f:ode17}
\end{figure}
{\bf 4.} It remains to prove the convergence $G^{t,\bar x}_{p,n}\to G^{t,\bar x}_p$
in $\L^1_{loc}(\R)$. 

As a preliminary, we remark that in the diffusion semigroups some trajectories
can move backward, from some $y_k$ to $y_{k-1}$.   However, the probability that 
this happens goes to zero as $n\to\infty$.  

Next, we observe that our previous assumption $y_i\leq\bar x<y_{i+1}$ implies the identities
\bel{Gpp}\bega{cl}
G^{t,\bar x}_p(y)\,=\,G^{t,\bar x}_{p-1}(y)\qquad &\hbox{for}\quad y\leq y_{j-1}\,,\\[3mm]
G^{t,\bar x}_{p,n}(y)\,=\,G^{t,\bar x}_p(y)~=~1\qquad&\hbox{for}\quad  y\geq y_{j}\,,\enda\eeq
where $j=i+p$.

To  prove the convergence in $\L^1\bigl([y_{j-1},y_j]\bigr)$, 
 for any point $z\in \,]y_{j-1}, y_j[\,$ and  any $\ve>0$,
we will show that 
\bel{Gnn}
 \limsup_{n\to \infty} ~G^{t,\bar x}_{p,n}(z-\ve)-\ve~\leq~
 G^{t,\bar x}_p(z)~\leq~\liminf_{n\to \infty} ~G^{t,\bar x}_{p,n}(z+\ve)+\ve.\eeq
 
 In analogy with (\ref{Gan}), we now set
\bel{Ga}\bega{rl}\Lambda _n(\tau)&\doteq~\mu_{1,n}^{\tau, y_{j-1}}\bigl([0,\tau]\times [z, y_j]\bigr),\\[3mm]\Lambda (\tau)&\doteq~\mu_1^{\tau, y_{j-1}}\bigl([0,\tau]\times [z, y_j]\bigr).
\enda\eeq
In other words, $\Lambda_n(\tau)$ is the probability that a trajectory, starting
at $y_{j-1}$, either stops at $y_j$ or reaches a point $\geq z$ at time $\tau$. 
 
We now have
\bel{Lam}
G^{t,\bar x}_{p}(z)~=~1-
\int_0^t \Lambda(t-\zeta) \, d F^{t,\bar{x}}_{p-1,j-1}(\zeta),\eeq
and an entirely similar formula holds for $G^{t,\bar x}_{p,n}(z)$.
For any $\ve>0$, the convergence $F^{t,\bar{x}}_{p-1,n,j-1}(\zeta)\to F^{t,\bar{x}}_{p-1,j-1}(\zeta)$
and $\Lambda_n\to \Lambda$ implies (\ref{Gnn}). This proves the desired $\L^1$ convergence
$G^{t,\bar x}_{p,n}\to G^{t,\bar x}_{p}$ in $\L^1([y_{j-1}, y_j]\bigr)$.
\v
{\bf 5.} Finally we observe that, for a fixed $t>0$, when $p$ is sufficiently large one has
$\tau^p(\bar x, \omega)>t$ with probability 1, for all $\bar x\in \R$.
Indeed, since trajectories of the semigroup are monotone, taking $p=N$ would suffice.
 Hence  $\mu^{t, \bar x}_N$ is supported on $\{t\}\times \R$, 
 and the weak convergence $\mu^{t,\bar x}_{N,n}\wto  \mu^{t,\bar x}_N$
 implies the weak convergence 
 (\ref{EP}).
 \endproof

\vs

\section{Approximating a Markov semigroup with smooth diffusions}
\setcounter{equation}{0}
\label{s:6}
In this section we complete the
proof of Theorem~\ref{t:2}. This  will require several  steps. 
\v
{\bf 1.} It will suffice to construct diffusions of the form
\bel{gen} dX_t~=~g_n(X_t)\, dt + \sigma_n\, dW_t\eeq
where $g_n$ is piecewise constant.  Indeed, when this is done,  for each $n\geq 1$
we can  then take a mollification:
$f_n= \phi_\delta *  g_n$, where the mollifier satisfies
$$\phi_\delta\in \C^\infty_c\,,\qquad \hbox{Supp}(\phi_\delta) = [-\delta, \delta], \qquad\qquad
\delta <\!< \sigma_n\,.$$
Without loss of generality we can assume that the
Markov semigroup satisfies the properties listed in Lemma~\ref{l:41}. Namely\begi
\item[{\bf (P)}] {\it The
 function $f$ is piecewise constant with compact support, the sets $\S,S^*, \Omega^*$ are finite and disjoint, and $\mu$ is the zero measure.  }
 \endi

Indeed, assume that  Theorem~\ref{t:2} holds true in this special case. 
To prove that the theorem remains valid for a general Markov semigroup, 
 let $(f_n)_{n\geq 1}$ be  a sequence of piecewise constant functions with compact support, constructed as
in Lemma~\ref{l:41}. For each $n\geq 1$,  there exists a sequence of diffusion processes of the form 
\[
dX~=~f^{(k)}_n(X)\, dt + \sigma^{(k)}_n dW,
\]
with $f^{(k)}_n\in \C^\infty$, such that $f^{(k)}_n\rightarrow f_n$ in the sense of graph
and  $\sigma^{(k)}_n\rightarrow 0$ as $k\to\infty$.  Moreover,
for every bounded continuous function $\vp\in \C^0(\R)$,
\[
\lim_{k\to\infty} 
\int \vp(x)P^{(n,k)}_t(x_0,dx)~=~\int \vp(x) P^{(n)}_t(x_0, dx).
\]
Consider the set of rational times $\mathcal{T}\doteq \mathbb{Q}\cap [0,T]$, 
and observe that the product set $\mathcal{T}\times\bigl(\mathbb{Q}\cup\S\cup\S^*\cup\Omega^*\bigr)$ is countable.
We can then select a sequence  $(k_n)_{n\geq 1}$ such that, setting 
$$\ov \sigma^{(n)}\,\doteq \,\sigma_n^{(k_n)},\qquad
\overline{P}^{(n)}_t\doteq P^{(n,k_n)}_t,$$  we have $\ov \sigma^{(n)}\to 0$.
Moreover, 
for every bounded continuous function $\vp\in \C^0(\R)$, we can achieve
\bel{conv1}
\lim_{n\to\infty} 
\int \vp(x)\overline{P}^{(n)}_t(x_0,dx)~=~\int \vp(x) P_t(x_0, dx)\qquad\forall (t,x_0)\in\mathcal{T}\times\bigl(\mathbb{Q}\cup\S\cup\S^*\cup\Omega^*\bigr).
\eeq
Equivalently, for every $(t,x_0)\in\mathcal{T}\times\bigl(\mathbb{Q}\cup\S\cup\S^*\cup\Omega^*\bigr)$, one has
\bel{distrib}
\lim_{n\to\infty}\overline{P}^{(n)}_t\bigl(x_0,\,]-\infty,x]\bigr)~=~P_t\bigl(x_0,\,]-\infty,x]\bigr)\qquad
\hbox{for a.e.~} x\in \R.
\eeq
We claim that the above convergence remains valid for all $(t, x_0)\in [0,T]\times\R$.
Indeed, for any  $t\in\mathcal{T}$,  the functions
 $$z_0~\mapsto~ P_t\bigl(z_0, \,]-\infty,x]\bigr),\qquad 
 z_0~\mapsto~ \overline{P}^{(n)}_t\bigl(z_0, \,]-\infty,x]\bigr)$$ are monotone increasing,  
and the map $z_0\mapsto P_t\bigl(z_0, \,]-\infty,\cdot ]\bigr)$ is continuous 
 w.r.t.~the $\L^1$ distance, at every $z_0\in \R\backslash (\S\cup\S^*\cup\Omega^*)$.
Using the fact that $\mathbb{Q}\cup\S\cup\S^*\cup\Omega^*$ is dense in $\R$, 
an approximation argument shows that the limit  (\ref{distrib}) remains valid 
for every $x_0\in\R$. Therefore, (\ref{conv1}) holds for every $(t,x_0)\in\mathcal{T}\times \R$. Finally, to obtain (\ref{conv1}) for every $t\in [0,T]$, we  observe that by {\bf (A1)}
the speed of all trajectories is bounded by $M$. Hence
\bel{conv3}\sup_{\|\vp\|_{\C^1}\leq 1} 
\int \vp(x)\Big[ P_t(x_0,dx)-P_s(x_0,dx)
\Big]~\leq~M\,|t-s|.
\eeq
\bel{conv4}\sup_{\|\vp\|_{\C^2}\leq 1} 
\int \vp(x)\Big[ \overline{P}^{(n)}_t(x_0,dx)-\overline{P}^{(n)}_s(x_0,dx)
\Big]~\leq~M\,|t-s| + M_1 \cdot \ov \sigma^{(n)},
\eeq
for a suitable constant $M_1$.
This establishes the convergence (\ref{EP}) for every $(t,x_0)\in [0,T]\times \R$
and every $\vp\in \C^2(\R)$. By an approximation argument, the same holds for every 
$\vp\in \C^0(\R)$.
\v
{\bf 2.}  It now remains to give the proof, in the case where the 
Markov semigroup satisfies the properties  {\bf (P)}.
As in (\ref{x0N}), we call $x_0<\cdots<x_N$ the points in
$\S\cup S^*\cup \Omega^*$, together with all points where
$f$ has a jump. Thanks to Lemma \ref{St2}, it suffices to show that, after inserting the stopping times
at all points $y_j$ in (\ref{y1N}), the corresponding transition kernels converge, for every 
initial point $\bar x$ and every $t>0$.
We are thus left with the task  of constructing a sequence of diffusion approximations
(\ref{gen}) on an interval $J=[y_{j-1}, y_j]$ where the function $f$ is piecewise constant, with a single jump at the interior point $x_j$.    In the easy 
case where $x_j\notin \S\cup S^*\cup \Omega^*$,  it suffices to choose $g_n(x) = f(x)$ for every $x\in J$ and $n\geq 1$.  The three main remaining cases will be discussed in the following steps.
\v
{\bf 3.}   CASE 1: $x_j\in \S$.  
Without loss of generality, we assume  $x_j=0$.  
Also, to fix ideas, assume that  the restriction of $f$ to the interval $J$ 
is given by (\ref{bode}), with $a,b>0$. See Fig.~\ref{f:ode18}, left.
Other cases can be treated similarly.

To make sure that the trajectories of the vanishing diffusion limit 
stop at the origin, for $x\in J$  we then define
\bel{gbx}
g_n(x)~=~\left\{  \bega{cl}  a\quad &\hbox{if}\quad  x< -\sqrt{\sigma_n}\\[2mm]
 0\quad &\hbox{if}\quad |x|\leq \sqrt{\sigma_n}\,,\\[2mm]
b\quad &\hbox{if}\quad x> \sqrt{\sigma_n}\,.\enda\right.\eeq
To prove that this sequence of diffusions achieves  the desired limit  (\ref{EP}),  let $\tau_n$ be the first time when a random trajectory of (\ref{gen}) starting at
the origin  reaches one of the two points  $\pm \sqrt{\sigma_n}$.  
For every given $t>0$, one has
\bel{t-0}\bega{rl}
\P\bigl\{\tau_n\leq t\bigr\}&\ds\leq~
\P\left(\sup_{s\in [0,t]} |W(s)|~\geq~{1\over \sqrt{\sigma_n}}\right)\\[4mm]
&\ds \leq~4\cdot \left(1-{1\over \sqrt{2\pi t}}\int^{{1\over\sqrt{\sigma_n}}}_{0}e^{{-x^2\over 2t}}dx\right)~\xrightarrow{\sigma_n\to 0}~0.
\enda\eeq
In the case where the motion starts at $x_0=0$, 
this yields 
\[
\lim_{n\to\infty}P^{(n)}_t\Big(  0\,;~  [-\sqrt{\sigma_n}\,,~\sqrt{\sigma_n}]\Big)~=~1.
\]
Therefore (\ref{EP}) holds for $x_0=0$. 
\medskip

Next, assume  $x_0\in J$ but $ x_0>0$.   
In this case, the result can be achieved by standard results in large deviation theory.
Since we are assuming $f(x)=b>0$ for $x>0$,
on the interval $[\sqrt \sigma_n, y_j]$ 
the diffusion process takes the form
 \bel{difb}dX_t ~=~b \, dt + \sigma_n \, dW_t\,,\qquad\quad X_0=x_0\,.\eeq
 
For any $\delta\in \,]0,x_0]$ and $T>0$, as $\sigma_n\to 0$ the probability that a random trajectory $X_t(\omega)$ starting at 
$x_0$ reaches a distance $>\delta$ from the limit solution $x(t)= x_0 + bt$ satisfies
$$\P\Big\{|X_t(\omega)- x(t)|~>\delta\quad\hbox{for some}~t\in [0,T]\Big\}
~\xrightarrow{\sigma_n\to 0}~0.$$
This already implies (\ref{EP}).
\medskip

Finally, consider the case where $x_0<0$.
For any given $0<\delta<|x_0|$, consider the time
$$T\,= \, {-\delta - x_0\over a}.$$
As before, the probability that   a random trajectory 
$X_t(\omega)$ starting at 
$x_0$ reaches a distance $>\delta$ from the limit solution $x(t)= x_0 + at$ satisfies
$$\P\Big\{|X_t(\omega)- x(t)|~>\delta\quad\hbox{for some}~t\in [0,T]\Big\}
~\xrightarrow{\sigma_n\to 0}~0.$$
Next, for a random trajectory such that 
$$\bigl|X_T(\omega)-x(T)\bigr|~=~\bigl|X_T(\omega)+\delta\bigr|~<~\delta,$$
let $\tau_n>T$ be the first time where $X_t(\omega)$ touches a point in the set
$\bigl\{-3\delta, \delta\bigr\}$.    As soon as $\sqrt{\sigma_n}<\delta$, the same 
argument used at (\ref{t-0}) yields
\bel{t-1}
\P\bigl\{\tau_n\leq T+t\bigr\}~\xrightarrow{\sigma_n\to 0}~0.
\eeq
Since $\delta>0$ was arbitrary, this again yields the weak convergence (\ref{EP}).
\v
{\bf 4.} 
CASE 2: 
 $x_j\, \in \, \Omega^*$ is a point from which both an increasing and a decreasing trajectory 
 initiate.   Again, we assume that $x_j=0$, and  that
 the restriction of $f$ to the interval $J$ is given by 
 (\ref{bode}), with $a<0<b$.   
Moreover, let $\theta\in [0,1]$ be the probability that a random trajectory
starting at the origin moves to the right.
 For each given $\sigma_n>0$, we consider the diffusion process (\ref{gen}), 
where
\bel{gn*}
g_n(x)~=~\left\{
\bega{rl} a\quad &\hbox{if} ~~x< \xi_n\,,\\[2mm]
b\quad &\hbox{if} ~~x> \xi_n\,,\enda\right.\eeq
for a suitable point  $\xi_n$, with $\xi_n\to 0$ as $n\to\infty$.   

To achieve the weak convergence (\ref{EP}) for the  initial point $x_0=0$,  
the points $\xi_n$ must be carefully chosen.
Toward this goal,
consider again the piecewise constant function $f$ in (\ref{bode}).
Performing the change of variables $y= {x/ \sqrt{\sigma_n}}$,  
the diffusion process
\bel{diff*}dX_t~=~f(X_t)\, dt + \sigma_n\, dW_t\eeq
becomes
\bel{diff}
d\Tilde{X}_t~=~\Tilde{f}_{n}(\Tilde X_t)dt+\sqrt{\sigma_n}\, dW_t\,,
\eeq
where the new drift is
\bel{gnn*}
\Tilde f_n(y)~=~\begin{cases}
{a/ \sqrt{\sigma_n}}~~~&\mathrm{if}\qquad  y<0,\\[2mm]
{b/ \sqrt{\sigma_n}}~~~&\mathrm{if}\qquad  y>0.
\end{cases}
\eeq
Let $\Tilde{\tau}^y_n$ be the first time when a random trajectory of (\ref{diff}) starting from 
 a point $y\in [-1,1]$ hits one of the two points in the set $\{-1,1\}$. 
Consider the expected values
\bel{v-u-tid}
\left\{\bega{rl}
\Tilde{v}_n(y)&\doteq ~E\left(\Tilde{\tau}^y_n\right)\\[2mm]
\Tilde{u}_n(y)&\doteq~E\left(\Tilde{X}_{\Tilde{\tau}^y_n}\right).\enda\right.
\eeq
It is well known that
 the functions $\Tilde{v}_n,\Tilde{u}_n$ provide the unique solutions to the boundary value problems
\bel{bvv} \left\{ \bega{rll}\Tilde f_n(y) v_y + \ds{\sigma_n\over 2} v_{yy}&=\,-1 \qquad &\hbox{if}~~ y\in \,]-1,1[\,,\\[3mm]
v&=~0\qquad &\hbox{if}~~  y\in \{-1,1\},\enda\right.\eeq
\bel{bvu} \left\{ \bega{rll}\Tilde f_n(y) u_y +\ds 
{\sigma_n\over 2} u_{yy}&=~0 \qquad &\hbox{if}~~  y\in \,]-1,1[\,,\\[3mm]
u(y)&=~y\qquad &\hbox{if}~~ y\in \{-1,1\}.\enda\right.\eeq
For a proof, see for example Chapter 9 in  \cite{O}.

The solutions of  (\ref{bvv}) and (\ref{bvu}) have the form
\bel{Tv}
\Tilde{v}_n(y)~=~\begin{cases}
-\ds {(y+1)\cdot \sqrt{\sigma_n}\over a}+c_{1,n}\cdot \left(e^{-a_{n}y}-e^{a_{n}}\right)&\qquad \mathrm{if}\quad -1<y\leq 0,\\[3mm]
\ds {(1-y)\cdot \sqrt{\sigma_n}\over b}+c_{2,n}\cdot \left(e^{-b_{n}y}-e^{-b_{n}}\right)&\qquad \mathrm{if}\quad 0\leq y<1,
\end{cases}
\eeq
\bel{Tu}
\Tilde{u}_n(y)~=~\begin{cases}
-1+2\cdot {\ds b\left(e^{a_{n}}-e^{-a_{n}y}\right)\over \ds b(e^{a_n}-1)+a (1- e^{-b_{n}})}&\qquad \mathrm{if}\quad -1<y\leq 0,\\[4mm]
1+\ds 2\cdot{a\left(e^{-b_{n}}-e^{-b_{n}y}\right)\over b(e^{a_n}-1)+
a (1- e^{-b_{n}})}&\qquad \mathrm{if}\quad 0\leq y<1,
\end{cases}
\eeq
and 
where the constants $a_n, b_n, c_{1,n}, c_{2,n}$ satisfy
\bel{apc}
a_n~\doteq~{2a\over \sigma_n^{3/2}},\qquad b_n~\doteq~{2b\over \sigma_n^{3/2}}\qquad\mathrm{and}\qquad |c_{1,n}|,|c_{2,n}|~\leq~{C\over \sqrt{\sigma_n}}.
\eeq
We observe that $\Tilde{u}_n$ is strictly increasing in $[-1,1]$. 
Moreover, as $n\to +\infty$ we have the limits
\bel{lsi} \sup\,\left\{ \Tilde{u}_n(x)\,;~x\in \Big[-1, \, -{1\over 3} \sqrt{\sigma_n}\Big]\right\}~\to ~-1,
\qquad \inf\,\left\{ \Tilde{u}_n(x)\,;~x\in \Big[{1\over 3} \sqrt{\sigma_n}\,,~1\Big]\right\}~\to ~1.\eeq
Therefore, for  any given  $0<\theta<1$ and every $n\geq 1$ large enough, we can uniquely choose an intermediate point $\zeta_n\in \big[-\sqrt{\sigma_n},\sqrt{\sigma_n}\big]$ such that 
\bel{x-n}
\Tilde{u}_n(\zeta_n)~=~ 2\theta-1.
\eeq
Moreover, from (\ref{apc}) and (\ref{Tv}) it follows
\bel{time0}
\lim_{n\to\infty} \Tilde{v}_n(x)~=~0\qquad\forall x\in [-1,1].
\eeq
Indeed, since the drifts in (\ref{diff})-(\ref{gnn*}) become very large, the average time needed 
for a random trajectory to exit from the interval $[-1,1]$ approaches zero.

Going back to the original space variable $x$,  we now set 
$\xi_n = \sqrt{\sigma_n}\, \zeta_n$.  According to the previous analysis,
for a random trajectory of (\ref{gen}), (\ref{gn*}),
starting at the origin,
the following holds.
\begi
\item[(i)] The random time $\tau_n(\omega)$ at which this trajectory reaches one of the points
$\xi_n \pm \sqrt{\sigma_n}$ approaches zero as $\sigma_n\downarrow 0$.
\item[(ii)] The probabilities of reaching the right or the left point are
given by
$$\bega{rl}\P\Big\{ X_{\tau_n}(\omega)\,=\,\xi_n+ \sqrt{\sigma_n}\Big\}&=~\theta,\\[2mm]
\P\Big\{ X_{\tau_n}(\omega)\,=\, \xi_n - \sqrt{\sigma_n}\Big\}&=~1-\theta.\enda
$$
\endi
Thanks to (\ref{lsi}), for all $n$ large enough 
we can assume $|\xi_n|\leq \ds{1\over 3} \sqrt{\sigma_n}$.    
A computation entirely similar to (\ref{t-0}) now shows that, as soon as a random trajectory
has reached one of the points $\xi_n \pm \sqrt{\sigma_n}$, the probability that
it goes back to the point $\xi_n$,  approaches zero as $\sigma_n\downarrow 0$.

To achieve the proof of weak convergence of the transition kernels starting at the origin, 
let $\delta>0$ be given.  Consider any random trajectory $X_t(\omega)$ of (\ref{gen}), with $g_n$ defined in (\ref{gn*}),
that starts at the point $\xi_n+\sqrt{\sigma_n}$ and then never touches the point $\xi_n$.
As an easy consequence of the theory of large deviations, for any $T>0$ we have
$$\P\left\{ \sup_{t\in [0,T]}~\bigr|X_t(\omega)-\xi_n+\sqrt{\sigma_n} + bt\bigr|~>~\delta
\Big|~X_0(\omega)= \xi_n+\sqrt{\sigma_n}\right\}~\xrightarrow{\sigma_n\to 0}~0.$$
Combining this fact with the above properties (i)-(ii), the weak convergence of the transitions kernels starting from $x_0=0$ is proved.

For every other initial point $x_0\in [y_{j-1}, y_j]\setminus \{x_0\}$, the weak convergence of the transition kernels is trivial.

\v
{\bf 5.}  CASE 3: $x_j\, \in \, \S^*$ is a point where trajectories of the Markov semigroup
stop for a random time  $T(\omega)\geq 0$
with Poisson distribution
$$\P\{ T(\omega)>s\}~=~e^{-\lambda s},$$
and then start moving again. As usual, we assume that $x_j=0$, while the drift $f$ is 
piecewise constant, as in (\ref{bode}).   To fix ideas, let $a,b>0$.   
The case $a,b<0$ is entirely similar.

We claim that this process can be approximated by a sequence of diffusions  as in 
(\ref{gen}), where each $g_n$ has the form
 \bel{gndef}g_n(x)~=~\left\{\bega{cl} 
 a\qquad&\hbox{if}
 \quad x<0,\\[2mm]
 -\eta_n\qquad&\hbox{if}
 \quad x\in [0, \ve_n],\\[2mm]
b \qquad&\hbox{if}\quad  x>\ve_n,\enda\right.\eeq
for a suitable choice of the constants $\ve_n,\eta_n$, with 
\bel{seq}0~<~\ve_n~<\!<~\sigma_n~<\!<~\eta_n.\eeq

Introducing the rescaled space variable 
 $y=\ds{ x\over \ve_n}$, the equation (\ref{gen})  becomes
\bel{dps}
dY_t~=~\tilde g_n(Y_t)dt+\tilde{\sigma}_n dW_t,
\eeq
with
\bel{tes}
\tilde g_n(y)~=~\left\{\bega{cl} 
\ds {a/ \ve_n}\qquad&\hbox{if}
 \quad y<0,\\[2mm]
\ds -\tilde{\eta}_n\qquad&\hbox{if}
 \quad y\in [0, 1],\\[2mm]
\ds {b/\ve_n} \qquad&\hbox{if}\quad y>1,\enda\right.
\qquad
\qquad (\tilde{\eta}_n,\tilde{\sigma}_n)~=~\left({\eta_n\over \ve_n},{\sigma_n\over \ve_n}\right).
\eeq

\begin{figure}[ht]
\centerline{\hbox{\includegraphics[width=8cm]{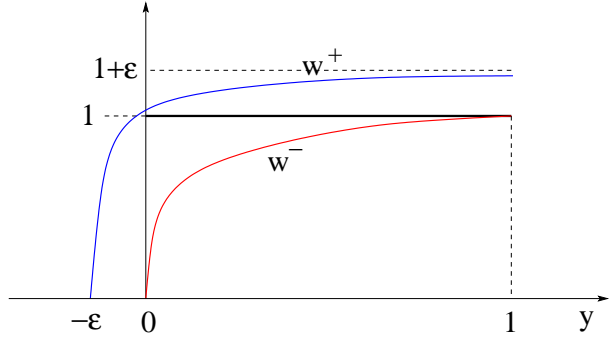}}}
\caption{\small  The  function $w$ at (\ref{w-}), used to construct a lower solution to (\ref{ze1}), and 
the functions 
$w_\delta$ at  (\ref{wdel})-(\ref{supsol}), used to construct upper solutions.
}
\label{f:ode19}
\end{figure}

{\bf 6.} Toward the analysis of (\ref{dps})-(\ref{tes}), in this step we perform a preliminary computation.
Let  $Y_t$ be the random solution to  
the diffusion process with constant coefficients, on the unit interval:
\bel{dse}
dY_t~=~-\eta \,Y_t dt+\sigma\, dW_t\eeq
 starting at the origin,
with the following boundary conditions. 
\begin{itemize}
\item {\bf Reflecting at the point $y=0$:} when $Y_t(\omega)= 0$, the particle is
reflected back inside the domain $]0,1[$;
\item  {\bf Absorbing at the point $y=1$: } when $Y_t(\omega)= 1$, its motion stops forever.
\end{itemize}
It is well known that in this case the distribution function
$$u(t,y)~=~\P\{ Y_t(\omega)\leq y\}$$
satisfies the parabolic equation
\bel{peq}
u_t~=~\eta u_y+ {\sigma^2\over 2} u_{yy}\,,
\eeq
with initial and boundary conditions
\bel{ibc3} u(0,y)~=~1,\qquad \qquad
\left\{ 
\bega{rl} \ds u(t,0)&=~0,\\[3mm]
u_x(t,1)&=~0.\enda\right.
\eeq
We seek a lower and an upper solution of (\ref{peq})-(\ref{ibc3}) 
in terms of suitable eigenfunctions. This leads us to
the boundary value problem
\bel{eig1}
{\sigma^2\over 2} w''(y) + \eta w'(y)~=\, -\lambda w(y),\qquad w(0)=0,\quad w'(1)=0\,.\eeq
Explicit solutions of  (\ref{eig1}) are found by solving
\[
{\sigma^2\over 2} r^2+\eta r+\lambda~=~0,\qquad\qquad r~=~{-\eta\pm\sqrt{\eta^2-2\sigma^2\lambda}\over \sigma^2}
~\doteq~ \gamma\pm s.
\]
For a given $\lambda>0$, choose any  $s>0$ and let $\eta,\sigma>0$, $\gamma<0$ be the constants such that 
\bel{beta-s}
{\sigma^2\over 2}~=~{(e^{2s}-1)^2\over 4s^2e^{2s}}\cdot\lambda,\qquad \eta~=~s\cdot {e^{s}+e^{-s}\over e^{s}-e^{-s}}\cdot\sigma^2\qquad\mathrm{and}\qquad \gamma~=~\ds{-\eta\over \sigma^2}~=~-s\cdot {e^{2s}+1\over e^{2s}-1}.
\eeq
In this case, an increasing solution to (\ref{eig1}), normalized so that $w(1)=1$, is explicitly found to be (see Fig.~\ref{f:ode19})
\bel{w-}
w(y)~=~{e^{(\gamma+s)y}-e^{(\gamma -s)y}\over e^{\gamma+s} -e^{\gamma-s}}\qquad\forall y\in [0,1].
\eeq
The function
\bel{lsol}
w^-(t,y)~=~e^{-\lambda t} w(y)\eeq
thus provides a lower solution to the parabolic problem (\ref{peq})-(\ref{ibc3}).

\begin{figure}[ht]
\centerline{\hbox{\includegraphics[width=9cm]{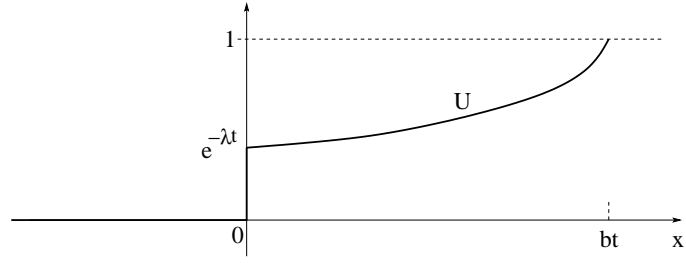}}}
\caption{\small The transition kernel (\ref{tker}), corresponding to a random Poisson waiting time at $x=0$,
followed by motion with constant speed $b>0$.}
\label{f:ode20}
\end{figure}

We now work toward the construction of an upper solution.
We claim that, for  every $\delta>0$ small, one can 
a function $w=w_\delta(x)$ such that (see Fig.~\ref{f:ode19})
\bel{wdel}w(-\delta)=0,\qquad w(0)= 1,
\qquad w(1)\leq 1+\delta, \qquad w'(1)>0,\eeq
and moreover
\bel{supsol}
{\sigma^2\over 2} w'' +\eta w' + (\lambda-\delta) w~\leq ~0\qquad
x\in [-\delta,1].\eeq
To construct the upper solution $w_\delta$,
let $\lambda,s>0$ be given, and let $\sigma, \eta,\gamma$ be as in (\ref{beta-s}).
For any small $\delta>0$, define the constants $s_{\delta},\lambda_{\delta}>0$  
in terms of the two equations 
\bel{s-ve}
\gamma~=~-s_{\delta}\cdot {e^{2s_{\delta}(1+2\delta)}+1\over e^{2s_{\delta}(1+2\delta)}-1}\qquad\mathrm{and}\qquad \lambda_{\delta}~=~{\gamma^2-s_{\delta}^2\over 2\sigma^2}.
\eeq
Then, on the larger interval $[-\delta,1]$,  the  problem
\bel{eig3}
{\sigma^2\over 2} w''(y) + \eta w'(y)~=\, -\lambda_{\delta} w(y),\qquad w(-\delta)=0,\qquad w(0)~=~1,\eeq
has a unique increasing solution   (see Fig.~\ref{f:ode19}), namely
\bel{wep}
w_{\delta}(y)~=~  {e^{(\gamma+s_{\delta})(y+\delta)}-e^{(\gamma-s_{\delta})(y+\delta)}\over e^{(\gamma+s_{\delta})\delta}-e^{(\gamma-s_{\delta})\delta} }\,, \qquad\qquad  y\in [-\delta,1].
\eeq
Using (\ref{s-ve}), we compute
\[
w_{\delta}(1)~=~{e^{(\gamma+s_{\delta})(1+\delta)}-e^{(\gamma-s_{\delta})(1+\delta)}\over e^{(\gamma+s_{\delta})\delta}-e^{(\gamma-s_{\delta})\delta} }~\leq~{e^{\gamma+s_{\delta}}\over 1-e^{-2s_{\delta}\delta}},
\]
\bel{w'(1)>0}
w'_{\delta}(1)~=~{e^{(\gamma-s_\delta)(1+\delta)}\over e^{(\gamma+s_{\delta})\delta}-e^{(\gamma-s_{\delta})\delta}}\cdot (\gamma-s_\delta)\cdot \left(e^{-2s_{\delta}\delta}-1\right)~>~0.
\eeq
Since $\ds s\cdot {e^{2s}+1\over e^{2s}-1}=s_{\delta}\cdot {e^{2s_{\delta}(1+2\delta)}+1\over e^{2s_{\delta}(1+2\delta)}-1}$, one has
\[
s\cdot {e^{2s}+1\over e^{2s}-1}\leq s_{\delta}\cdot {e^{2s_{\delta}}+1\over e^{2s_{\delta}}-1},\qquad s_{\delta}~\leq~s\cdot \left(1+{2\over e^{2s}-1}\right)~\leq~s+1.
\]
In particular, for $s\geq 1$ and $0< \delta\leq 1/2$, we have 
\[
s~\leq~s_{\delta}~\leq~2s,\qquad \gamma+s_{\delta}~=~-{2s_{\delta}\over e^{2s_{\delta}(1+2\delta)}-1}~\leq~-2se^{-8s},
\]
and hence
$$
w_{\delta}(1)~\leq~{e^{-2se^{-8s}}\over 1-e^{-2s\delta}},\qquad \qquad {-3s^4e^{2s}\over  (e^{2s}-1)^2\lambda}~\leq~{s^2-s^2_{\delta}\over 2\sigma^2}~=~\lambda_{\delta}-\lambda~\leq~0.
$$
Therefore, by choosing $s=\delta^{-2}$, so that $\delta=\ds{1/ \sqrt{s}}$, one obtains
\bel{w(1)<1+ve}
w_{\delta}(1)~\leq~1+\delta,\qquad  \qquad\lambda-\delta~\leq~\lambda_{\delta}~\leq~\lambda,
\eeq
provided that  $s>0$ is sufficiently large.
\v
{\bf 7.}  We are now ready to construct a sequence of  
diffusion processes with piecewise constant drift, which approximate a Poisson 
waiting time at  $x_0=0$.   For
each $t>0$, the transition kernel we need to approximate is (see Fig.~\ref{f:ode20})
\bel{tker}U(t,x)~=~\prob\{ X_t(\omega)\leq x\}~=~\left\{ \bega{cl} 0\qquad &\hbox{if}\quad x\leq 0,\\[2mm]
\ds \exp\left\{ -\lambda\Big( t-{x\over b}\Big)\right\}  \qquad &\hbox{if}\quad x \in [0,  bt], \qquad  \\[2mm]
1\qquad &\hbox{if}\quad x \geq bt.\enda\right.\eeq
Indeed, since particles travel with constant speed $b>0$, 
a particle will reach the point $x>0$ after 
time $t$ if and only if it departs from the origin after time $t-(x/b)$. The probability of this event
is $\exp\left\{ -\lambda\bigl( t-{x\over b}\bigr)\right\}$.

Next, we compare the kernel (\ref{tker})  with the solution
 $u=u_n(t,x)$ of  the parabolic problem 
\bel{ze1}u_t+ g_n(x) u_x~=~{\sigma_n^2\over 2} u_{xx} \,,\qquad x\in \R,
\eeq 
with $g_n$ as in (\ref{gndef}),  and with initial data
\bel{ze2}
u(0,x)~=~\left\{ 
\bega{rl} 0\quad \hbox{if} ~~x<0\,,\cr
1\quad \hbox{if}~~x\geq 0\,.\enda\right.\eeq
We claim that, for a suitable choice of the sequences, $\sigma_n, \ve_n,\eta_n\to 0$, one can achieve 
the weak convergence of the transition kernels.  In other words, the sequence of solutions
$u_n=u_n(t,\cdot)$ of the parabolic problems (\ref{ze1})-(\ref{ze2}) converge to $U(t,\cdot)$ in $\L^1_{loc}(\R)$, for every $t>0$. 
The claim will be proved by constructing suitable sequences of upper and lower solutions.

The sequences $\ve_n, \eta_n\to 0$, are chosen as follows.
First, we take  a sequence $s_n\to +\infty$.  Then, in view of (\ref{tes}) and (\ref{beta-s}), we define
\bel{setan}
{\Tilde \sigma_n^2\over 2}~=~{(e^{2s_n-1})^2\over 4s_n^2e^{2s_n}}\cdot\lambda,\qquad \Tilde \eta_n~=~s_n\cdot {e^{s_n}+e^{-s_n}\over e^{s_n}-e^{-s_n}}\cdot \Tilde \sigma^2_n\,,\eeq
\bel{etsin}\eta_n~=~\ve_n \Tilde\eta_n,\qquad\qquad \sigma_n~=~\ve_n \Tilde\sigma_n\,.
\eeq
Notice that here
$\Tilde \sigma_n, \Tilde\eta_n \to +\infty$. 
Therefore, we need to choose a sequence $\ve_n\to 0$ converging  to zero fast enough 
so that  $\eta_n, \sigma_n\to 0$.

\v

\begin{figure}[ht]
\centerline{\hbox{\includegraphics[width=9cm]{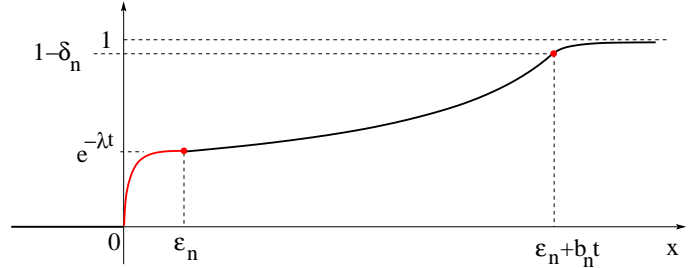}}}
\caption{\small  The lower solution defined at (\ref{lowso}).}
\label{f:ode21}
\end{figure}
Let $w_n^-$ be the solution to (\ref{eig1}) with $\sigma=\tilde{\sigma}_n$, $\eta=\tilde{\eta}_n$ and $w_n^-(1)=1$. 
We now set $\delta_n= 1/\sqrt{s_n}$.   For suitable sequences
$\kappa_n\to +\infty$ and
$b_n\downarrow b$, we define
\bel{lowso}u^-_n(t,x)~=~\left\{\bega{cl} 0\qquad &\hbox{if}\quad x<0,\\[3mm]
(1-\delta_n) e^{-\lambda t} w_n^{-}(x/\ve_n)\qquad &\hbox{if}\quad x\in [0,\ve_n],\\[3mm]
\ds (1-\delta_n) \exp\left\{ -\lambda\Big( t-{x-\ve_n\over b_n}\Big)\right\}  \qquad &\hbox{if}\quad x \in 
\bigl[\ve_n, \ve_n+ b_n t\bigr], \qquad  \\[3mm] 1- \delta_n \exp\Big\{-\kappa_n \bigl(x-\ve_n -b_n t\bigr)\Big\}
\qquad &\hbox{if}\quad x \geq \ve_n+ b_n t.\enda\right.\eeq
This  lower solution is illustrated in Fig.~\ref{f:ode21}.  
Here we choose $\kappa_n\to +\infty$ fast enough so that $\delta_n\kappa_n\to +\infty$.
Finally, we choose  the decreasing sequence $b_n\downarrow b$  
converging slowly enough
so that the traveling profile
$$v(t,x)~=~1-\delta_n e^{-\kappa_n(x-\ve_n-b_n t)}$$
is a lower solution. This will be the case if  
$$v_t~\leq~{\sigma_n^2\over 2 } v_{xx} + b_n v_x\,.$$
That is, iff
$$- b_n \kappa_n ~\leq~-{\sigma_n^2\over 2} \kappa_n^2 - b \kappa_n\,.$$
We thus need
$$b_n-b~\geq~{\sigma_n^2\kappa_n\over 2}\,.$$
Notice that, since the product $\delta_n\kappa_n\to +\infty$, this ensures that
the derivative of the solution at the matching point $x=\ve_n+b_n t$ jumps upward.
The fact that this derivative also jumps upward at $x=0$ and at $x=\ve_n$ is obvious.

\begin{figure}[ht]
\centerline{\hbox{\includegraphics[width=10cm]{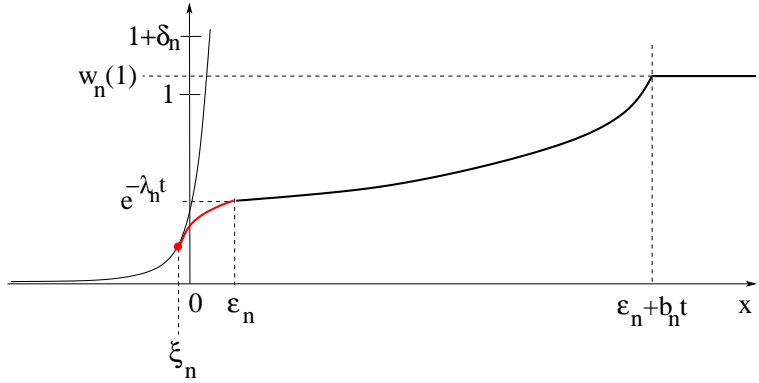}}}
\caption{\small  The upper solution defined at (\ref{upso}).}
\label{f:ode22}
\end{figure}

Next, an upper solution is constructed as follows. As before, set $\delta_n=1/\sqrt{s_n}$.
Let $w_n = w_{\delta_n}$ 
be a sequence of functions satisfying (\ref{wdel})-(\ref{supsol}).

We then consider the rescaled functions 
$ x\mapsto w_n(\ve_n^{-1} x)$, 
choosing a sequence $\ve_n\downarrow 0$  fast enough so that  $\ve_n^{-1} w'_n(1)\to +\infty$.  

Next, we take a sequence $\kappa_n\to +\infty$. For each $n\geq 1$,  we choose the values $x_n$ 
and  $\xi_n\in [-\delta_n \ve_n, 0]$
so that at $x=\xi_n$ the two functions
$$x\mapsto w_n(\ve_n^{-1} x)\qquad\hbox{and}\qquad x\mapsto {e^{\kappa_n (x-x_n)}}$$
are tangent:
$$w_n(\ve_n^{-1} \xi_n)~=~ e^{\kappa_n (\xi_n-x_n)},\qquad \ve_n^{-1} w_n'(\ve_n^{-1} \xi_n)~=~ \kappa_n e^{\kappa_n (\xi_n-x_n)}\,.$$
Finally, we choose an increasing sequence of speeds $b^-_n\uparrow b$ and define
(see Fig.~\ref{f:ode22})
\bel{upso}u^+_n(t,x)~=~\left\{\bega{cl} e^{-\lambda_n t} e^{\kappa_n (x-x_n)}
\qquad &\hbox{if}\quad x<\xi_n,\\[3mm]
e^{-\lambda_n t} w_n(x/\ve_n)\qquad &\hbox{if}\quad x\in [\xi_n,\ve_n],\\[3mm]
\ds w_n(1) \cdot \exp\left\{ -\lambda_n\Big( t-{x-\ve_n\over b^-_n}\Big)\right\}  \qquad &\hbox{if}\quad x \in 
\bigl[\ve_n, ~\ve_n +  b^-_n t\bigr], \qquad  \\[3mm]
w_n(1)\qquad &\hbox{if}\quad x \geq \ve_n +  b^-_n t.\enda\right.\eeq

The condition $\ve_n^{-1} w_n(1)\to +\infty$ guarantees that at the point $x=\ve_n$ the derivative $\partial_x u^+_n$ has a downward jump.    On the other hand, this derivative is continuous at the matching point 
$x=\xi_n$, and has a downward jump at $x= \ve_n + b^-_nt$.    
We thus conclude that $u^+_n$ is an upper solution.

Since the sequences $u^-_n$ and $u^+_n$ at (\ref{lowso}) and (\ref{upso}) 
both converge to the distribution function $U$ at (\ref{tker}), for all $t\geq 0$,
this  proves the weak convergence of the transition kernels, starting at the origin.

For a starting point  $x_0\not=0$, the convergence of the corresponding transition kernels
is a straightforward.
\v
{\bf 8.} To cover the set of all possibities, one should also consider
a dynamics of the form (\ref{bode}) in two additional cases:

CASE 4: $a,b>0$, and there is no stopping time at $x=0$.

CASE 5: $b<0<a$, so that all trajectories stop forever, as they reach the origin.

In both of these cases, the approximation of the Markov semigroup
with a diffusion process is trivial. Indeed, taking  $g_n(x) = f(x)$
for all $n\geq 1$,    as the coefficients $\sigma_n\downarrow 0$,
the transition kernels for the diffusion process converge to the corresponding transition kernels of the Markov semigroup.

In view of Lemma~\ref{St2}, this completes the proof of Theorem~\ref{t:2}.
\endproof

\end{document}